\documentstyle[txmac,a4,
amssymb,%
case,%
twoside,%
nocaphead,%
epsf,%
mypic,times,mathptm]{article}

\advance\oddsidemargin by -1.8cm
\advance\evensidemargin by -1.8cm
\advance\textwidth by 3.6cm

\def\mynewtheo#1#2{%
\newtheorem{@#1}{#2}
\newenvironment{#1}{\begin{@#1}\rm}{\end{@#1}}}

\mynewtheo{lemma}{Lemma}
\mynewtheo{exer}{Exercise}
\mynewtheo{theo}{Theorem}
\mynewtheo{rem}{Remark}
\mynewtheo{defi}{Definition}
\mynewtheo{conj}{Conjecture}
\mynewtheo{corr}{Corollary}
\mynewtheo{proposition}{Proposition}
\mynewtheo{question}{Question}
\mynewtheo{exam}{Example}

\newenvironment{theorem}{\begin{theo}}{\end{theo}}
\newenvironment{conjecture}{\begin{conj}}{\end{conj}}

\newenvironment{eqn}{\begin{equation}}{\end{equation}}

\begin{document}

\makeatletter

\def\bysame{\same[\kern2cm]\,}
\def\lfra{\leftrightarrow}
\def\qed{\hfill\@mt{\Box}}
\def\@mt#1{\ifmmode#1\else$#1$\fi}
\def\qqed{\hfill\@mt{\Box\enspace\Box}}

\let\ap\alpha
\let\tl\tilde
\let\sg\sigma
\let\dl\delta
\let\Dl\Delta
\let\eps\varepsilon
\let\gm\gamma
\let\bt\beta
\def\cF{{\cal F}}
\def\bZ{{\Bbb Z}}
\def\bN{{\Bbb N}}
\def\spn{\mbox{\operator@font span}\,}
\let\ds\displaystyle
\def\cf{\text{\rm cf}\,}
\def\md{\min\deg}
\def\Md{\max\deg}

\def\epsfs#1#2{{\catcode`\_=11\relax\ifautoepsf\unitxsize#1\relax\else
\epsfxsize#1\relax\fi\epsffile{#2.eps}}}
\def\epsfsv#1#2{{\vcbox{\epsfs{#1}{#2}}}}
\def\vcbox#1{\setbox\@tempboxa=\hbox{#1}\parbox{\wd\@tempboxa}{\box
  \@tempboxa}}
\def\lz{\linebreak[0]\verb}

\def\@test#1#2#3#4{%
  \let\@tempa\go@
  \@tempdima#1\relax\@tempdimb#3\@tempdima\relax\@tempdima#4\unitxsize\relax
  \ifdim \@tempdimb>\z@\relax
    \ifdim \@tempdimb<#2%
      \def\@tempa{\@test{#1}{#2}}%
    \fi
  \fi
  \@tempa
}

\def\go@#1\@end{}
\newdimen\unitxsize
\newif\ifautoepsf\autoepsftrue

\unitxsize4cm\relax
\def\epsfsize#1#2{\epsfxsize\relax\ifautoepsf
  {\@test{#1}{#2}{0.1 }{4   }
		{0.2 }{3   }
		{0.3 }{2   }
		{0.4 }{1.7 }
		{0.5 }{1.5 }
		{0.6 }{1.4 }
		{0.7 }{1.3 }
		{0.8 }{1.2 }
		{0.9 }{1.1 }
		{1.1 }{1.  }
		{1.2 }{0.9 }
		{1.4 }{0.8 }
		{1.6 }{0.75}
		{2.  }{0.7 }
		{2.25}{0.6 }
		{3   }{0.55}
		{5   }{0.5 }
		{10  }{0.33}
		{-1  }{0.25}\@end
		\ea}\ea\epsfxsize\the\@tempdima\relax
		\fi
		}

\let\old@tl\~
\def\~{\raisebox{-0.8ex}{\tt\old@tl{}}}

\author{A. Stoimenow\footnotemark[1]\\[2mm]
\small Max-Planck-Institut f\"ur Mathematik,\\
\small Vivatsgasse 7, D-53111 Bonn, Germany,\\
\small Postal Address: P.\ O.\ Box: 7280, D-53072 Bonn,\\
\small e-mail: {\tt alex@mpim-bonn.mpg.de},\\
\small WWW: {\hbox{\tt http://guests.mpim-bonn.mpg.de/alex}}
}

{\def\thefootnote{\fnsymbol{footnote}}
\footnotetext[1]{Supported by a DFG postdoc grant.}
}

\title{\large\bf \uppercase{On the crossing number of positive knots
and braids and}\\[2mm]
\uppercase{braid index criteria of Jones and
Morton-Williams-Franks}\\[4mm]
{\small\it This is a preprint. I would be grateful
for any comments and corrections!}}

\date{\large Current version: \curv\ \ \ First version:
\makedate{2}{10}{2000}}

\maketitle

\long\def\@makecaption#1#2{%
   \vskip 10pt
   {\let\label\@gobble
   \let\ignorespaces\@empty
   \xdef\@tempt{#2}%
   }%
   \ea\@ifempty\ea{\@tempt}{%
   \setbox\@tempboxa\hbox{%
      \fignr#1#2}%
      }{%
   \setbox\@tempboxa\hbox{%
      {\fignr#1:}\capt\ #2}%
      }%
   \ifdim \wd\@tempboxa >\captionwidth {%
      \rightskip=\@captionmargin\leftskip=\@captionmargin
      \unhbox\@tempboxa\par}%
   \else
      \hbox to\captionwidth{\hfil\box\@tempboxa\hfil}%
   \fi}%
\def\fignr{\small\sffamily\bfseries}%
\def\capt{\small\sffamily}%

\newdimen\@captionmargin\@captionmargin2cm\relax
\newdimen\captionwidth\captionwidth\hsize\relax

\let\reference\ref
\def\eqref#1{(\protect\ref{#1})}

\def\proof{\@ifnextchar[{\@proof}{\@proof[\unskip]}}
\def\@proof[#1]{\noindent{\bf Proof #1.}\enspace}

\def\myfrac#1#2{\raisebox{0.2em}{\small$#1$}\!/\!\raisebox{-0.2em}{\small$#2$}}
\def\abstractname{}

\@addtoreset {footnote}{page}

\renewcommand{\section}{%
   \@startsection
         {section}{1}{\z@}{-1.5ex \@plus -1ex \@minus -.2ex}%
               {1ex \@plus.2ex}{\large\bf}%
}
\renewcommand{\@seccntformat}[1]{\csname the#1\endcsname .
\quad}

{\let\@noitemerr\relax
\vskip-2.7em\kern0pt\begin{abstract}
\noindent{\bf Abstract.}\enspace
We give examples of knots with some unusual properties of
the crossing number of positive diagrams or strand number of
positive braid representations. In particular we show that
positive braid knots may not have positive minimal (strand number)
braid representations, giving a counterpart to results of
Franks-Williams and Murasugi. Other examples answer questions of
Cromwell on homogeneous and (partially) of Adams on almost
alternating knots.
\\[2mm]
We give a counterexample to, and a corrected version of a theorem of
Jones on the Alexander polynomial of 4-braid knots. We also give an
example of a knot on which all previously applied braid index criteria
fail to estimate sharply (from below) the braid index. A relation
between (generalizations of) such examples and a conjecture of Jones
that a minimal braid representation has unique writhe is discussed.
\\[2mm]
Finally, we give a counterexample to Morton's conjecture
relating the genus and degree of the skein polynomial.
\end{abstract}
}

{
\parskip 1pt plus 1pt minus 1pt\relax
\tableofcontents}

\section{Introduction}

The braid index $b(L)$ of a knot or link $L$ is defined to be
the minimal number of strands of a braid, whose closure is the link.
(That such a braid always exists was first shown by Alexander
\cite{Alexander}.) To determine the braid index of $L$, one is
seeking general lower and upper estimates on $b(L)$. Upper estimates
can be obtained in the obvious way by writing down braid
representations of $L$ (although finding a braid representation
realizing $b(L)$ may be sometimes difficult), so the harder problem
is to estimate $b(L)$ from below. Very little was known in general on
this problem (except some early results of Murasugi \cite{Murasugi3}
for 3-braids, which required much effort), until Jones
discovered his polynomial invariant in 1984 \cite{Jones}.
His construction made heavy use of braid
representations, and thus he obtained several conditions for
knots, in particular of low braid index. Briefly later the HOMFLY
(skein) polynomial \cite{HOMFLY,LickMil} was discovered, which gave rise
to the presently most commonly used braid index estimate\footnote{In
the sequel we will be interested only in lower estimates.} for
$b(L)$, the MWF inequality \cite{Morton,WilFr}. This inequality
determined the braid index of
all knots in \cite{Rolfsen}, except five. They were dealt with by the
3-braid formula of Murakami \cite[corollary 10.5]{Murakami}, or by
applying the MWF inequality on their 2-cable \cite{MorSho2}.

One of the central points of this paper is the study
of these braid index inequalities with particular regard to
positive braids.

We discuss a conjecture of Jones in \cite{Jones2},
that a minimal braid representation has unique writhe, and relate
this conjecture to the MWF inequality and its cabled versions
\cite{MorSho2}. A consequence of this relation is that
on a counterexample to Jones's conjecture \em{any} cable version
of the MWF inequality will fail to estimate sharply the braid index.
Therefore, at the present state of the art it is very unlikely to
find (that is, to prove some link to be) a counterexample
to Jones's conjecture, except possibly if it is a 4-braid\footnote{%
It has been claimed by Birman \cite{BirMen} that the truth of this
conjecture for 4-braids follows from Jones's work \cite{Jones2}, but
this claim is possibly incorrect; see \S\reference{S7}.}.
On the quest of such an example we found knots for
which both the MWF inequality and its 2-cable version
(and hence any previously applied method) fail to estimate sharply
the braid index. We will show one of these knots.

Also, we consider one of Jones's original
criteria in \cite{Jones}. We provide a counterexample to it, showing
that it needs correction, and we give the corrected version. 

Then we turn to Jones's unity root criteria for the Jones
polynomial. We give an example showing that these criteria sometimes
can estimate the braid index better than the MWF inequality,
and thus deserve (although apparently neglected after MWF)
to be considered in their own right.

Another aim of the paper, which we will begin with, is to show some
examples of positive knots with with unusual behaviour of classical
invariants as braid index and crossing number in positive braid
representations and diagrams.

Positive knots are called the knots with diagrams of all crossings
positive (see e.g. \cite{MorCro}). This class of knots
contains as subclass the braid positive
knots, those which are closures of positive braids\footnote{Some
authors very confusingly call `positive knots' what we will
call here `braid positive knots'.}. Such knots
were studied in knot theory, \em{inter alia} because of their
relevance to the theory of singularities \cite{Hirzebruch}
and dynamical systems \cite{BirWil}. Thus they received much
attention in previous publications \cite{Busk,Cromwell2,Rudolph}.

Positive and braid positive knots have been studied in many papers
jointly with alternating (braid) knots as
a subclass of the homogeneous knots and braids \cite{Cromwell,%
Stallings}. It is now known (see e.g. \cite{Murasugi,Murasugi2}) that
reduced alternating diagrams are of minimal crossing number
and that reduced alternating braid representations are of minimal
strand number. Simple examples show that neither of this is true
for positive/homogeneous diagrams/braid representations,
so that the reasonable question is whether there always exists
at least some such minimal diagram/braid representation.
A partial positive answer in the case of positive braids was given
in \cite{WilFr} for positive braids containing a full twist.
In \cite{Murasugi} it was remarked that a positive/homogeneous braid
representation of minimal strand number has also minimal crossing
number. (Thus the positive answer for braid representations for
a given link implies a positive answer for diagrams for this link.)
Here we show that the answers to both questions are in general negative.

\begin{theorem}
There exist knots with positive/homogeneous diagrams but with no
positive/homogeneous diagrams of minimal crossing number, or
with positive/homogeneous braid representations,
but with no positive/homogeneous braid representations of minimal 
strand number.
\end{theorem}

Beside these examples, we will prove some relations between the
crossing number and genus of braid positive knots. These inequalities
will enable us to show that certain knots, like Perko's knot $10_{161}$,
have no positive braid representations, or that the reduced positive
braid representation of some others, like the closed 4-braid
$(\sg_1\sg_3\sg_2^2)^3\sg_2$, is unique.

In the final section, we will give examples settling two conjectures
on possible inequalities between the genus and the degrees of the
skein polynomial, one of which is a 15-year-old problem of Morton
\cite{Morton4}.

Most of the examples presented below were found by examining the
tables of the program KnotScape of Hoste and Thistlethwaite
\cite{KnotScape}. Beside providing access to these tables, the program
offers the possibility to calculate their polynomial invariants
and to identify a knot in the table from a given diagram.
These features were used to large extent in the calculations
described below.

\subsection{Definitions and notation.} 

The $n$-strand braid group $B_n$ is generated by the elementary (Artin)
braids $\sg_i$ for $1\le i\le n-1$ with relations $\sg_i\sg_{i+1}\sg_i=
\sg_{i+1}\sg_i\sg_{i+1}$, called henceforth Yang-Baxter (YB)
relation, and $[\sg_i,\sg_j]=1$ for $|i-j|>1$ (the brackets denoting the
commutator), called commutativity relation.

By $\hat\ap$ we denote link, which is the braid closure of $\ap$.
Markov's theorem (see e.g. \cite{Morton3}) says that when
$\hat\ap=\hat\ap'$ then $\ap$ and $\ap'$ can be transformed into each
other by a sequence of conjugacies in the $B_n$'s
and moves of the type $\ap\lfra\ap\sg_{n}^{\pm1}\in B_{n+1}$ for
$\ap\in B_n$. We call the `$\to$' part of this move (which augments
the strand number by 1) stabilization, and its inverse destabilization.

By $[\bt]$ we denote the exponent sum of a braid $\bt$, that
is, the image of $\bt$ under the homomorphism $B_n\to\bZ$ given
by $[\sg_i]:=1$ for any $i$. A braid word $\bt$ is called positive
if its length equals $[\bt]$, and a braid $\bt$ is positive if it has a
positive word representation. (As already apparent, we will often abuse
the distinction between braids and braid words, as this will cause no
confusion.) By $[\bt]_i$ we denote the exponent sum
of the generator $\sg_i$ in the braid \em{word} $\bt$, which
is clearly not invariant under the YB relation. That is, $[\,.\,]_i$
is a homomorphism of the free group in the $\sg_j$ given by
$[\sg_j]_i:=\dl_{ij}$ (where $\dl$ is the Kronecker delta).

By $P_K$ or $P(K)$ we denote the
skein polynomial of $K$ \cite{HOMFLY}, and by $v$ its non-Alexander
variable. The span of $P$ in $v$ means the difference between its
maximal and minimal
degree in this variable. These degrees are denoted by $\Md_vP_K$ and
$\md_vP_K$, respectively. The other (Alexander) variable of $P$
is denoted by $z$. By $\max\cf_zP$ we mean the maximal coefficient of
$z$ in $P$ (which is a polynomial in $v$), i.e. the coefficient of
the maximal degree of $z$ in $P$.

The braid index $b(K)$ of a knot $K$ is defined by 
\[
b(K)=\min\{\,n\,|\, \exists\bt\in B_n\,:\,\hat\bt=K\,\}\,.
\]
A lower bound for the braid index
is given by the inequality of Franks--Williams \cite{WilFr} and Morton
\cite{Morton}: 
\[
b(K)\ge \mbox{$v$-span}\ P_K/2+1\,.
\]
The inequality of Morton--Williams--Franks will be subsequently
abbreviated as `MWF inequality' or simply as `MWF', and its
right hand-side will be called the `MWF bound' for $K$.

Whenever we talk of a diagram to be minimal, we always mean
minimality with respect to the crossing number of the knot it
represents. Similarly, a minimal braid representation is
meant with respect to the braid index of its closure (i.e.,
that the strand number of the braid realizes this braid index).

`W.l.o.g.' will abbreviate `without loss of generality' and
`r.h.s' (resp.\ `l.h.s') `right hand-side' (resp.\ `left hand-side').

A final remark on knot tables and notation is in order. It is understood
that alternative work on knot tabulation to that of Hoste,
Thistlethwaite and Weeks is being done by Aneziris
\cite{Aneziris}. Unfortunately, it seems like every new knot tabulator
chooses and insists on his own numbering convention for knots, which
will (and, in fact, already did\footnote{I know of at least
two cases where the Hoste--Thistlethwaite numbering
of the Rolfsen knots was the origin of embarrassing confusions in
published material.}) lead to confusion in using the
different knot tables. It appears
most correct to stick to the convention of the first tabulator for each
crossing number. We use here the convention of Rolfsen's tables
\cite{Rolfsen} for $\le 10$ crossing knots and that of \cite{KnotScape}
for $\ge 11$ crossing knots, which coincides with those of the first
tabulators for any crossing number except 11, where the initial
(complete) tables were compiled by Conway \cite{Conway}. We apologize for
not using his numbering. An excuse is that all calculations have been
performed by KnotScape, which yet does not provide a translator between
its notation and that of Conway. For uniformity reasons, we will
need to continue using this convention in subsequent papers, too.

\section{Some interesting diagrams of $\bf 11_{550}$}

\subsection{A positive knot with no positive minimal diagram}

An intuitive question on positive knots (whose affirmation
in the alternating case was one of the big achievements of the
Jones polynomial) is whether any positive knot has a positive
minimal (crossing number) diagram (see \cite{pos}). This was known
to be true in the case the positive knot is alternating \cite{%
Nakamura} or of genus at most two \cite{gen2}. The following,
surprisingly simple, (non-alternating genus three) example 
shows that this need not be true in general.

{\begin{figure}[htb]
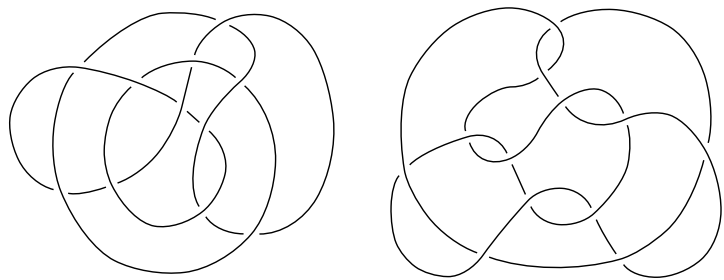

\[
\begin{array}{c@{\qquad}c}
\epsfsv{4.0cm}{t1-11-550} & \epsfsv{4cm}{t1-11-550-2}\\
\end{array}
\]
\caption{Two diagrams of the knot $11_{550}$. On the left its
(only) minimal diagram (which is not positive), and on the right a
positive 12 crossing diagram.\label{figpos}}
\end{figure}
}

\begin{exam}\label{11-550}
The knot $11_{550}$ has only one 11 crossing diagram shown on
the left of figure \reference{figpos}. The diagram is only
almost positive (i.e. has one negative crossing). However,
$11_{550}$ has a positive 12 crossing diagram shown on
the right. Thus it is positive, but has no positive minimal diagram.
As by \cite[theorem 4 and corollary 4.1]{Cromwell} any
homogeneous diagram of a positive knot must be positive,
this example simultaneously provides a negative answer to
question 2 in \cite[\S 5]{Cromwell}. (This also implies
a positive answer to question 1 therein, but this
answer was previously known to follow from the almost
positive diagram of the Perko knot~-- the mirror image of the
diagram of $10_{161}$ in the Rolfsen's tables \cite{Rolfsen}.)
\end{exam}

Another problem for positive knots is then in how far
the crossing number of a positive diagram can differ
from the crossing number of the knot. The presently
known result is obtained in \cite{pos} using the Gau\ss{}
sum theory of Polyak-Viro-Fiedler.

\begin{theorem}(\cite{pos})\label{thp}
If $D$ is a positive reduced diagram (i.e. with no nugatory
crossings) of a positive knot $K$ with $c(D)$ crossings,
then the crossing number $c(K)$ of $K$ satisfies $c(K)\ge\sqrt{2c(D)}$.
\end{theorem}

The bound is clearly not very sharp, and a much better
estimate appears to be true.

\begin{conjecture}
With the notation of theorem \reference{thp}, $c(K)\ge c(D)-2g(K)+1$,
where $g(K)$ is the genus of $K$.
\end{conjecture}

Here is some motivation for this conjecture.
\begin{itemize}
\item The conjecture is known to be true for $K$ fibered
\cite[corollary 5.1]{Cromwell}, of genus
at most $2$ \cite{gen2}, and for $c(D)\le 16$ by experiment
(note, that only the cases where $g(K)=3$ are relevant to check).

\item The inequality is sharp for all positive rational knots (i.e.,
is the best possible) \cite{fibo}, and also for some other knots,
e.g. $9_{16}$ (see figure \reference{fig9-16}). It is worth
remarking that all knots I found so far, for which the
inequality was sharp, are alternating and arborescent
(this was in particular always the case for $g(K)=3$ and
$c(D)\le 16$).
\end{itemize}

{\begin{figure}[htb]
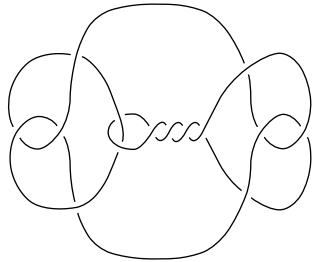

\[
\begin{array}{c}
\epsfsv{3.7cm}{t1-9-16-14}
\end{array}
\]
\caption{A positive 14 crossing diagram of the knot $9_{16}$.
\label{fig9-16}}
\end{figure}
}

A further question we can pose is

\begin{question}\label{q1}
Does a braid positive knot always have a (braid) positive
minimal diagram?
\end{question}

We will later answer this question positively for $\le 16$ crossing
knots, but also provide evidence against it by answering negative
question about closely related properties of braid positive knots.

\subsection{An example on almost alternating diagrams}

While discussing diagrams of the knot $11_{550}$, it is worth making an
aside from our positivity considerations to almost alternating diagrams.

Such diagrams were considered in \cite{Adams2} to be diagrams obtainable
from alternating diagrams by one crossing change, and almost
alternating knots are knots having such diagrams, but which are not
alternating. A surprising variety
of knots turns out to be almost alternating, in particular very
many low crossing number knots. Using Conway's description of
$\le 11$ crossing prime knots \cite{Conway} and a simple way to
manipulate their Conway notation, all such non-alternating knots
were found to be almost alternating with 3 exceptions. They are
shown on figure 5.54 of \cite{Adams}, and (from left to right
and in our notation) are $11_{550}$, $11_{485}$ and $11_{462}$.

A computer check showed that in fact our knot $11_{550}$ is
almost alternating. 2 almost alternating diagrams of it are shown
on figure \reference{figaal}. Such diagrams show that in
general it will be hard to decide on almost alternation of a given knot.
Neither strong obstructions are known, nor any effective
method for seeking almost alternating diagrams is available.
(It is not hard to see e.g. that many almost alternating knots
have infinitely many almost alternating diagrams, so that an upper
control on their crossing number is not possible.)

The other two knots are indeed problematic. We cannot prove them not
to be almost alternating, but no almost
alternating diagram was found after checking all
diagrams of $\le 16$ crossings and some diagrams of 17 crossings.

{\begin{figure}[htb]
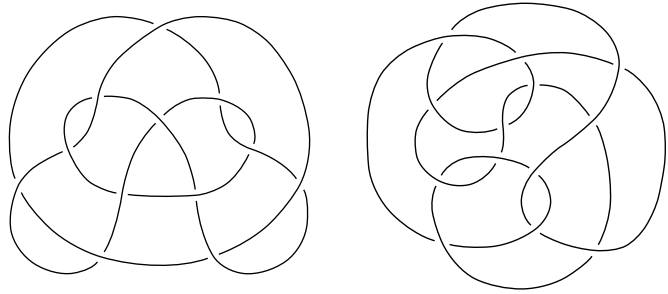

\[
\begin{array}{c@{\qquad}c}
\epsfsv{3.7cm}{t1-11-550-aal} & \epsfsv{4cm}{t1-11-550-aal2}\\
\end{array}
\]
\caption{Two almost alternating diagrams of the knot $11_{550}$.
\label{figaal}}
\end{figure}
}

\section{The crossing number and genus of positive braids}

\subsection{Irreducible positive braids}

In \cite[corollary 5.1]{Cromwell}
it was shown that for a positive $c$ crossing 
diagram of a fibered positive knot of genus $g$ it holds $g\ge c/4$.
This inequality in particular applies to braid positive
diagrams. We will improve the inequality by showing that
for braid positive knots $K$ we have $g(K)\ge c(K)/4+2$ with finitely
many exceptions. This is based on a continuation of the
investigation on the crossing number of irreducible positive
braids which was initiated in \cite{pos}.

\begin{defi}
For $n\ge 2$ define the number $d_n$ to be the minimal number of
crossings of a positive $n$-(strand )braid $\bt$ such that there
is no positive braid $\bt'$ of smaller crossing number, or
(because of the genus) equivalently strand number,
with $\hat\bt=\hat\bt'$, and $\hat\bt$ is a prime knot.
We call such a braid $\bt$ \em{irreducible}.
\end{defi}

\begin{rem}
In the definition no requirement is made (the crossing number $[\bt]$
of) $\bt$ to realize $c(\hat\bt)$. Anything that follows
would remain true with this modified definition, but the definition
as made above appears more natural because an answer to question
\reference{q1} is unclear.
\end{rem}

\begin{proposition}\label{pdn}
The values of $d_n$ for small $n$ are as follows.
\begin{eqn}\label{d_n}
\begin{array}{c*{5}{|c}}
n & 2 & 3 & 4 & 5 & 6 \\[0.3em]
\hline
\ry{1.2em}d_n & 3 & 8 & 11 & 16 & 19
\end{array}
\end{eqn}
\end{proposition}

We will later describe how they were obtained.

Here we show the following theorem.

\begin{theorem}\label{thdn}
\def\labelenumi{\theenumi)}\mbox{}\\[-18pt]
 
\begin{enumerate}
\item\label{item1} For $n\ge 7$ we have $d_n\ge 2n+6$.
\item\label{item2} $\ds d_n\le \frac{9}{2}n-
\left\{\begin{array}{cc}
3 & n\equiv 0\,(4) \\
\myfrac{9}{2} & n\equiv 1\,(4) \\
6 & n\equiv 2\,(4) \\
\myfrac{3}{2} & n\equiv 3\,(4)
\end{array}\right.\,=\,\frac{9}{2}n-\frac{3}{2}\bigl(\,(n+1)\bmod 4+1\,
\bigr)$\,.
\end{enumerate}
\end{theorem}

\newpage

\proof 
\begin{enumerate}
\def\labelenumi{\theenumi)}
\item
Take a braid $\bt$ realizing $d_n$ for $n$ fixed. For irreducibility
reasons we must have $[\bt]_i\ge 2$ for any $1\le i\le n-1$. Our
aim is to show that
\begin{eqn}\label{10}
\sum_{i=1}^3[\bt]_{n-i}\ge 10\,\mbox{\ \ and similarly\ \ }
\sum_{i=1}^3[\bt]_{i}\ge 10\,.
\end{eqn}
Then $[\bt]\ge 2n+6$ for $3<n-3$, i.e. $n\ge 7$.

When writing $\bt=\prod_j\sg_{i_j}$, we can modulo YB relations
assume that the index sum $\sum_ji_j$ is minimal. Using this word
representation for $\bt$, consider the
subword of $\bt$ made up of $\sg_{n-2}$ and $\sg_{n-1}$, keeping
separate parts separated by subwords of $\bt$ made up of $\sg_{n-i}$,
$i>2$. Thus we can write $\bt$ as
\[
\bigl(\sg_{n-2}^{a_{1,1}}\ap_{1,1}\sg_{n-2}^{a_{1,2}}\ap_{1,2}\dots \sg_{n-2}^{a_{1,n_1}}\ap_{1,n_1}\sg_{n-1}^{b_1}\bigr)\dots
\bigl(\sg_{n-2}^{a_{k,1}}\ap_{k,1}\sg_{n-2}^{a_{k,2}}\ap_{k,2}\dots \sg_{n-2}^{a_{k,n_k}}\ap_{k,n_k}\sg_{n-1}^{b_k}\bigr)\,,
\]
with the $\ap_{l,j}$ standing for subwords containing only $\sg_{n-i}$, $i>2$.

We can use commutativity relations to assure that each one of the
subwords $\ap_{l,j}$ contains at least one $\sg_{n-3}$, and that all
$n_i>0$. 

We have that $k\ge 2$, else $\bt$ decomposes. If $n_l=1$ and $a_{l,1}=1$
for some $l$, we can apply (after some commutativity relations)
a YB relation $\sg_{n-1}\sg_{n-2}\sg_{n-1}\to\sg_{n-2}\sg_{n-1}\sg_{n-2}
$ to reduce the index sum of the word, a contradiction to our
assumption. Thus assume that $\sum_{i=1}^{n_l}a_{l,i}\ge 2$ for
all $1\le l\le k$. Therefore,
\[
[\bt]_{n-2}\,=\,\sum_{l=1}^k\sum_{i=1}^{n_l}a_{l,i}\ge 4\,.
\]

\begin{caselist}
\case \label{csA}
Assume that $[\bt]_{n-2}=4$ and $[\bt]_{n-1}=2$. According
to the distribution of $\sg_{n-2}$ we have 3 possibilities.

\begin{caselist}
\case \label{cs1} $k=2$, $n_1=n_2=2$ and $a_{1,1}=a_{1,2}=a_{2,1}=a_{2,2}=1$.

\case \label{cs2} $k=2$, $n_1=2$, $n_2=1$, $a_{1,1}=a_{1,2}=1$, and
$a_{2,1}=2$.

\case \label{cs3} $k=2$, $n_1=n_2=1$ and $a_{1,1}=a_{2,1}=2$.
\end{caselist}

Case \reference{cs3} is excluded, because the closure is not connected
(i.e., not a knot).
In cases \reference{cs1} and \reference{cs2} the following argument
applies.

Since all $\ap_{i,j}$ contain the letter $\sg_{n-3}$, we have
\[
[\bt]_{n-3}\ge \sum_{i=1}^k\,n_i\,.
\]
Thus in case \reference{cs1} we have $[\bt]_{n-3}\ge 4$. If in
case \reference{cs2} $[\bt]_{n-3}=3$, then for one of $i=1,2$,
$\ap_{i,1}$ contains the letter $\sg_{n-3}$ exactly once.
Then (after some commutativity relations) the subword
$\sg_{n-2}^{a_{i,1}}\ap_{i,1}\sg_{n-2}^{a_{i,2}}$ can be made
to admit a YB relation 
$\sg_{n-2}\sg_{n-3}\sg_{n-2}\to\sg_{n-3}\sg_{n-2}\sg_{n-3}$,
a contradiction to our assumption. Thus $[\bt]_{n-3}\ge 4$,
and \eqref{10} holds. (The second inequality therein follows
analogously to the first one.)

\case \label{csB} Otherwise, $[\bt]_{n-1}+[\bt]_{n-2}\ge 7$, so
$\sum_{i=1}^3[\bt]_{n-i}\ge 9$. Again one needs to check that
the case $\sum_{i=1}^3[\bt]_{n-i}=9$ cannot occur.
For this one applies the same type of argument, but
the case list becomes too large to be effectively handled
manually, so one safer checks the cases by computer (see \S
\reference{cd} for more details on this calculation). 
\end{caselist}

%
%

\item
We write down explicit positive braids of the given number of
crossings. To show that they are irreducible, we use the value
$\big|\,V\left(e^{\pi i/3}\right)\,\big|$ (where $V$ is the
invariant introduced in \cite{Jones}) on their closure
and apply proposition 14.6 of \cite{Jones2}. It is easy to see that
this value is preserved by a $3$-move, which in the context
of braid words means cancelling subwords of the type $\sg_i^3$.
As all the braids we will write down become trivial after
a sequence of such cancellations, their closure satisfies
$\big|\,V\left(e^{\pi i/3}\right)\,\big|=\sqrt{3}^{n-1}$,
and thus the braids are irreducible (a more special type of this
argument was given in corollary 15.5 of \cite{Jones2}). To
show primeness, we use the result of \cite{Cromwell2}.

\begin{caselist}
\case $n\equiv 0\,(4)$. Consider 
\[
\left(\sg_1^3\sg_2\dots\sg_{n-3}^3\sg_{n-2}\sg_{n-1}^3\right)^2
\sg_2\sg_4\sg_6\dots\sg_{n-2}\,.
\]

\case $n\equiv 1\,(4)$. Consider 
\[
\bt=\left(\sg_1^3\sg_2\dots\sg_{n-2}^3\sg_{n-1}\right)^2
\sg_2\sg_4\sg_6\dots\sg_{n-1}\,.
\]

\case $n\equiv 2\,(4)$. Consider 
\[
\left(\sg_1\sg_2^3\dots\sg_{n-2}^3\sg_{n-1}\right)^2
\sg_1\sg_3\sg_5\dots\sg_{n-1}\,.
\]
\case $n\equiv 3\,(4)$. Consider
\[
\left(\sg_1\sg_2^3\dots\sg_{n-3}^3\sg_{n-2}\sg_{n-1}^3\right)^2
\sg_1^4\sg_3\sg_5\dots\sg_{n-2}\,.
\]
\end{caselist}

It is easy to check that in all cases the closures are knots and that
the crossing numbers are as stated above. \qed
\end{enumerate}

The following straightforward consequence shows in how far the
inequality $g\ge c/4$ for braid positive knots can be improved.

\begin{corr}
For almost all (i.e., all but finitely many) braid positive knots,
$g\ge c/4+2$. Moreover, for any constant $C$ with $g\ge C\cdot c$
for almost all braid positive knots we have $C\le \myfrac{7}{18}$. \qed
\end{corr}

\begin{rem}\label{rm1}
Although we have not yet proved the values for $d_n$ in \eqref{d_n},
the proof of \reference{item1}) already shows that $d_n\ge 2n+2$ for
$n=5,6$.
\end{rem}

\begin{exam}\label{161}
Theorem \reference{thdn} can be used to show that some fibered
positive knots, like Perko's, are not braid positive. $K=10_{161}$
has genus 3 and braid index $3$. If $K=\hat\bt$ with $\bt\in B_n$
positive and w.l.o.g.\ irreducible, then $[\bt]\ge \max(c(K),d_n)$, and 
\[
g(K)=\frac{[\bt]-n+1}{2}\ge \frac{\max(c(K),d_n)-n+1}{2}\,.
\]
However, by remark \reference{rm1}, the r.h.s. is at least $4$ for
any $n\ge 3$, a contradiction. We will later see that this simple
reasoning (given the numbers $d_n$ computed for enough small $n$),
does not always work.
\end{exam}

\proof[of proposition \reference{pdn}]
First one generates all braid positive $\le 16$ crossing
diagrams from the tables of \cite{KnotScape} and identifies braid
positive knots from them. This led to the values of $d_n$ for
$n\le 5$. For $n=6$ one needed to exclude 17 crossing braids.
This was done by generating a superset of all irreducible 17
crossing braids (see \S\reference{cd} for more details how this
was done) and identifying their closures. All knots had
$\le 16$ crossings. That there is a 19 crossing irreducible 6-braid
will be shown by example later (see example \reference{exbt}). \qed

\begin{rem}
It would be interesting how the sequence of $d_n$ continues.
Sloane \cite{Sloane} reports on two sequences starting as in
\eqref{d_n}. One is related to \cite{CSV} and the other one made 
up of numbers $\equiv 0,3\,(8)$. It would be surprising if the
answer were that simple, though.
\end{rem}

\subsection{Examples}

\subsubsection{Minimal crossing number positive braid diagrams}

Although it is desirable to push further results of the
above type, there are many difficulties in controlling
positive braid representations. We illustrate this by a series
of examples. 

\begin{exam}\label{exbt}
The inequality of \reference{item1} of theorem \reference{thdn}
cannot be improved by trivial means. As noted in \cite{pos},
the braids $\{\beta_n\,|\,n \text{ odd}\,\}$ with
\[
\beta_n=\bigl((\sg_1\sg_3\dots\sg_{n-4}\sg_{n-2}^3)\,
(\sg_2^3\sg_4\dots\sg_{n-3}\sg_{n-1})\bigr)^2
\]
admit no YB relation. A computer check shows that $\bt_5$
is indeed irreducible (one of the examples showing $d_5=16$),
while $\bt_7$ reduces to the 19 crossing 6-braid
\[
\bt'_7=\sg_1\sg_2^3\sg_4^2\sg_3\sg_5\sg_4^3\sg_1\sg_3
\sg_5\sg_2^3\sg_3\sg_4\,,
\]
thus showing $d_6=19$,
and that $\bt_9$ also reduces to a 6-braid, this time of 21 crossings.

Calculation of the Fiedler polynomial $\Dl$ \cite{Fiedler} 
(see also \cite{Morton2}, but \em{not} \cite{Alexander})
of any of the stabilizations $\bt''_7$ of $\bt'_7$ or 
cyclic permutations or flips ($\sg_i\lfra\sg_{n-i}$) thereof, we had
(writing just the honest polynomial part and using the notation
of \cite{LickMil}) $\Dl(\bt''_7)=[0]\ 2\ 0\ 3\ 0 \ 5$, while $\Dl(
\bt_7)=[0]\ 2\ 0\ 2\ 0\ 6$, so $\bt_7$ and $\bt''_7$ are not conjugate.
It is not clear whether $\bt_7$ is not conjugate to a
stabilization of another positive conjugate of $\bt'_7$. (There is
an algorithm to list up all such positive conjugates
\cite{MortonElrifai}, but it is too complex, especially as I have
no computer version of it.) Anyway, this example
provides some evidence against an easy version of Markov's theorem
for positive braids. In any case one cannot obtain every minimal
positive braid representation from a given one just by YB moves,
cyclic permutations and destabilization.
\end{exam}

The following example shows that the argument in example 
\reference{161} very often fails.

\begin{exam}\label{exx}
The two mutant knots on figure \reference{figgc}, $15_{203432}$ and
$15_{203777}$, are positive fibred of genus 6, as many other
15 crossing knots, which are braid positive. The diagrams on the figure
are the only 15 crossing diagrams of these knots, and they are not
braid diagrams. In fact it turns out that both knots are not
braid positive. To show this, it suffices to show that they have
no positive braid representations of 17 crossings and 6 strands,
and of 16 crossings and 5 strands. Indeed the knots did not occur
in such representations (see \S\reference{cd}).
This shows that just comparing genus and
crossing number as in example \reference{161} will in general not
suffice to exclude braid positivity. Many more (several hundred)
such examples (again of genus 6) occurred at crossing number 16.
\end{exam}

We will come back to this example a little later when we consider
mutation.


{\begin{figure}[htb]
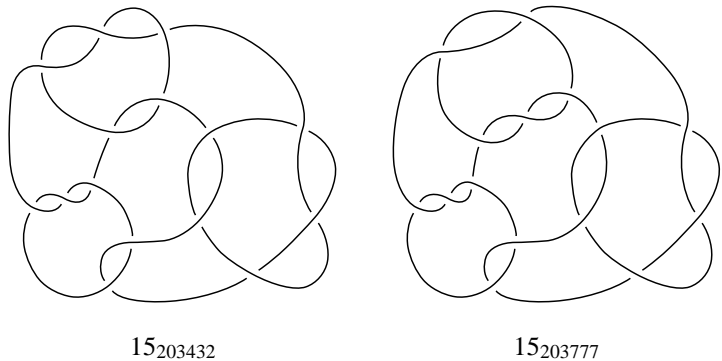

\[
\begin{array}{c@{\qquad}c}
\epsfsv{4cm}{t1-15-203432} & \epsfsv{4cm}{t1-15-203777}\\[19mm]
\ry{1.5em}15_{203432} & 15_{203777}
\end{array}
\]
\caption{Two knots which are fibered and positive, and have the
same genus and crossing number as other braid positive knots, but
are not braid positive.\label{figgc}}
\end{figure}
}

\begin{theorem}\label{pdo}
Any braid positive knot of $\le 16$ crossings has a braid positive
minimal diagram.
\end{theorem}

\proof We proceed as follows.

First we consider crossing number $\le 15$.
\begin{enumerate}
\def\labelenumi{\theenumi)}
\item\label{it1} We identified all knots with braid positive diagrams of $\le 15$
crossings, and found that all they have a braid positive
minimal diagram.

\item
Now we need to consider the cases, where the braid positive
knot has $\le 15$ crossings, but all its braid positive diagrams
have $>15$ crossings. We have from \cite{Cromwell,MurPrz}
that for a braid positive knot $K$, $\Md_zP_K=\md_vP_K(=2g(K))$
and $\max\cf_vP_K=z^{\md_vP_K}$. Thus we select all knots whose
$P$ polynomial has this property. If $K$ is braid positive
and $g(K)(=\Md_zP_K/2)\le 5$, then from theorem \reference{thdn}
and \eqref{d_n} we have a positive braid representation of $\le 13$
crossings, and checked that such knots have minimal braid
positive diagrams.

\item
If $g>5$, then the knot was among those identified in \reference{it1}),
except the knots in example \reference{exx}.
\end{enumerate}

Now consider crossing number 16.

Again we select the 16 crossing non-alternating knots satisfying
the above condition on $P$ with $g\ge 6$ ($g\le 5$ is dealt with
as above). If such knots are braid positive, then they must have
braid positive diagrams of $\le 17$ crossings (see example
\reference{exx}). There were 393 knots of $16$ crossings
with the $P$ condition for $g\ge 6$, which did not have
braid positive 16 crossing diagrams. All they had $g=6$
(i.e. no one had $g=7$), and no one of the 17 crossing 
braid positive diagrams (of genus 6) identified to any of these knots.
Thus they are all not braid positive. The same exclusion applied for the
two knots in example \reference{exx}. \qed

A final remark on minimal crossing number diagrams is that it is not
true that \em{any} minimal diagram of a braid positive knot is
braid positive. 

\begin{exam}
The knot $11_{444}$ is braid positive, but not all its minimal
diagrams are positive braid diagrams. (This phenomenon does
not occur for the braid positive knots in Rolfsen's tables.)
See figure \reference{11-444}.
\end{exam}

{\begin{figure}[htb]
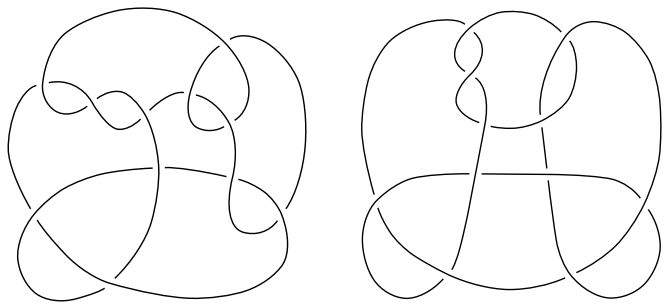

\[
\begin{array}{c@{\qquad}c}
\epsfsv{4.0cm}{t1-11-444-1} & \epsfsv{4cm}{t1-11-444-2}\\
\end{array}
\]
\caption{Two minimal diagrams of the knot $11_{444}$. On the left
a diagram, which is the diagram of a closed braid,
and on the right another diagram which is not. \label{11-444}}
\end{figure}
}

\subsubsection{Minimal strand number positive braid representations}

There is an observation of \cite{Murasugi} relating the crossing number
and braid index of braid positive knots. Murasugi showed a minimal
positive (or homogeneous) braid representation to have minimal crossing
number of its closure. He also showed that reduced alternating
braid representations are minimal. This is trivially not true for
positive braid representations (e.g. $(\sg_1\sg_2)^2$ for the trefoil),
so a natural question is whether there is at least one positive minimal
braid representation. In \cite{WilFr}, it was shown that
if $\bt$ factors as $\Dl^2\alpha$ with $\alpha$ positive and
$\Dl^2$ being the full twist braid (generating the center of the
braid group), then MWF is sharp for $\bt$, so that $\bt$
is minimal (this contains the case of torus knots considered in
\cite[proposition 7.5]{Murasugi}).

However, the following examples show that both conclusions are in
general problematic. Neither the MWF inequality can always help
to prove some positive braid representation of a given knot
to be minimal, nor needs such a minimal positive braid representation 
to exist at all.

\begin{exam}
Among braid positive 15 crossing knots, which (except the $(2,15)$-%
torus knot) have genus $6$ (and a positive 4-braid representation),
there are two knots with MWF bound $3$, see figure \reference{figMWF}.
That both knots are not closures of 3-braids can be shown using their
$Q$ polynomial \cite{BLM,Ho} and the formula of Murakami \cite{Murakami}
(see also \cite[theorem 2]{Kanenobu}). The (only) possible value for
the exponent sum of the hypothetic 3-braid representing the knot
can be found from the inequalities of \cite{Morton} and the
degrees of $P$ in the non-Alexander variable. Thus MWF need not
be sharp even on minimal positive braid representations.
(Another such example will be given later in \S\reference{SU}.)
\end{exam}

{\begin{figure}[htb]
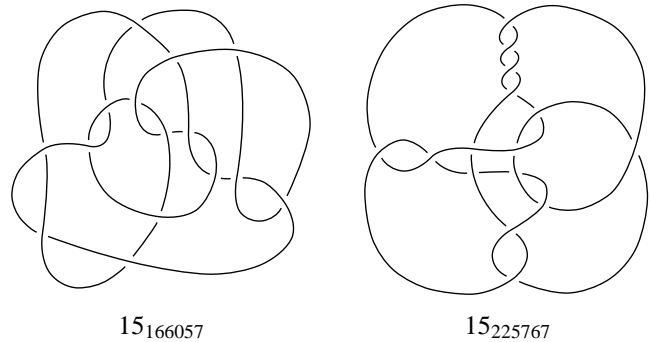

\[
\begin{array}{c@{\qquad}c}
\epsfsv{4cm}{t1-15-166057} & \epsfsv{3.8cm}{t1-15-225767}\\[3mm]
\ry{1.5em}15_{166057} & 15_{225767}
\end{array}
\]
\caption{Two braid positive knots with unsharp MWF.\label{figMWF}}
\end{figure}
}

\begin{exam}\label{exmbr}
Considering braid positive 16 crossing knots (they are of genus 6 or 7),
there are again two genus 6 knots (closures of positive 16 crossing
5-braids), $16_{472381}$ and $16_{1223549}$, whose MWF bound $4$
is unsharp on these braid representations (see figure
\reference{figmbr}). This time is was not
clear that both knots are not of braid index 4. The evidence against
this was strenghtened by their 2-cable-$P$, calculated with the
program of Morton and Short (see \cite{MorSho2}), whose $v$-span
was $14$ for both knots. Finally, an extensive check found indeed
$4$-braid representations of these knots, for example
\[
\sg_3^{-1}\sg_2\sg_3^3\sg_2^2\sg_1^2\sg_2^2\sg_1\sg_3\sg_2^2\sg_1^2\mbox{ and }
\sg_3^{-1}\sg_2\sg_3^2\sg_2^2\sg_1^2\sg_2^2\sg_1\sg_3\sg_2^2\sg_1^2\sg_2\,.
\]
Thus these two knots are braid positive, but have no positive minimal
braid representation. (By the same remark on Cromwell's work as in
example \reference{11-550}, these examples serve equally well also
for homogeneous braid representations.)
\end{exam}

{\begin{figure}[htb]
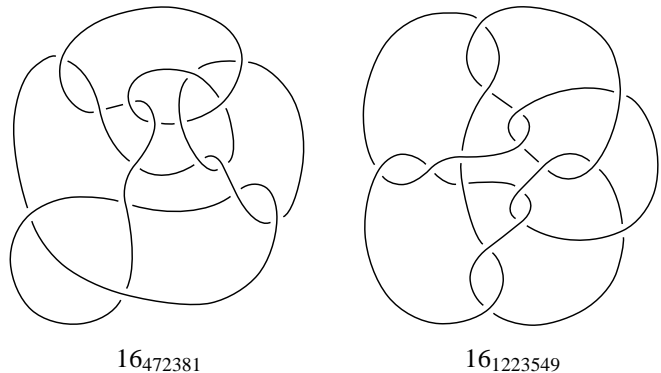

\[
\begin{array}{c@{\qquad}c}
\epsfsv{4cm}{t1-16-472381} & \epsfsv{4cm}{t1-16-1223549}\\[3mm]
\ry{1.5em}16_{472381} & 16_{1223549}
\end{array}
\]
\caption{Two braid positive knots with no minimal positive
braid representation.\label{figmbr}}
\end{figure}
}

Thus the question whether a braid positive braid index $n$
knot has a positive minimal braid representation has a negative
answer for $n=4$, and taking connected sums of these knots with
trefoils, also for $n>4$ (for this one needs to use e.g. the
result of \cite{Cromwell2}). The only non-settled case remains
$n=3$ (for $n=2$ the answer is positive by elementary means).

\begin{question}
If a knot $K$ has braid index $3$ and is the closure of a positive
braid, is it the closure of a positive $3$-braid?
\end{question}

Unfortunately, we know from the work of Birman \cite{Birman}
and Murakami \cite{Murakami}, that considering 3-braid knots
just via their polynomials $P$ (and hence $V$ and $\Dl$
\cite{Alexander}) and $Q$ will not
suffice to give a positive answer to this question.

In a similar way, taking iterated connected sums of the
knots in example \reference{exmbr} and using the result of 
Cromwell \cite{Cromwell2}, one obtains knots $K_n$ for which
the difference between $b(K_n)$ and the positive braid index
$b_p(K_n)$, the minimal strand number of a positive braid representation
of $K_n$, becomes arbitrarily large. It would be interesting to
find prime examples.

\begin{question}
Is there a sequence of prime braid positive knots $\{K_n\}$, for
which $b_p(K_n)-b(K_n)\to\infty$?
\end{question}

\subsubsection{Some examples on mutation}

The next examples concern mutation. Mutation was introduced by
Conway \cite{Conway}, and consists in building links which
turn out very difficult to distinguish. Mutation replaces in
a link diagram a tangle with its rotated version by $180^\circ$
around some axis (see \cite{LickMil}).
The following example shows that for braid positive knots, 
mutation may not be visible in positive braid diagrams.

{\begin{figure}[htb]
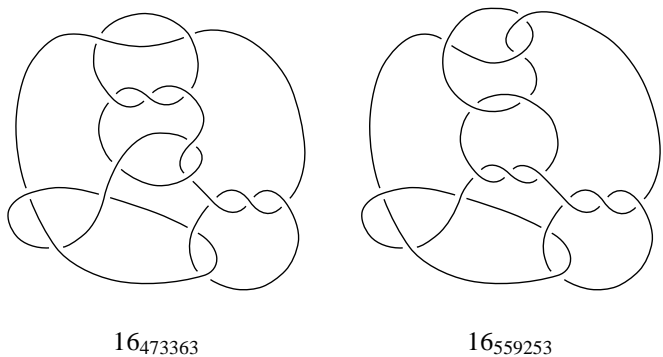

\[
\begin{array}{c@{\qquad}c}
\epsfsv{4cm}{t1-16_473363} & \epsfsv{4cm}{t1-16_559253}\\[19mm]
\ry{1.5em}16_{473363} & 16_{559253}
\end{array}
\]
\caption{Two braid positive mutants with no mutated positive braid
diagrams.\label{figmmu}}
\end{figure}
}

To fix a bit of terminology call a braid word \em{reduced}
if it has no isolated generator. (This is \em{not} to be confused
with the previous notion of irreducible, which is stronger.)

\begin{exam}
The knots $16_{473363}$ and $16_{559253}$ in figure \reference{figmmu}
are iterated mutants. It turns out that both knots have a unique 16
crossing diagram which is a diagram of a closed positive braid
(on 5 strands). The braids are
\begin{eqn}\label{iop}
\sg_1\sg_2^2\sg_3^3\sg_1^2\sg_4\sg_2^3\sg_3^3\sg_4
\mbox{ and }
\sg_1\sg_2^3\sg_3^3\sg_1\sg_4\sg_2^3\sg_3^3\sg_4\,.
\end{eqn}

We have already shown in the proof of
part \reference{item1} of theorem \reference{thdn}
(and its computational part to follow in \S\reference{xx}),
that for positive braids $\bt$ of genus 6
and $n\ge 7$ strands can be turned into such with isolated generators
by YB relations. If $\bt$ is itself reduced, then the YB relation
giving an isolated generator is of the form $\sg_1\sg_2\sg_1\to
\sg_2\sg_1\sg_2$ (or its version with $i$ replaced by $n-i$, $i=1,2$).
Thus after reducing the isolated generator $\sg_1$ (or $\sg_{n-1}$)
we obtain a braid with an edge generator ($\sg_1$ or $\sg_{n-1}$)
occurring in a square.

Moreover, it turns out that none of these knots occurs in a list
of 17 crossing 6-strand diagrams, which cannot be
reduced in the above way (they all represent 14 and 15
crossing knots; see \S\reference{ST3}). 

Thus, starting from a reduced braid representation of
these knots of $\ge 6$ strands, by iteratedly applying YB
relations and removing isolated (edge) generators, one must arrive
to a 5-strand representation with an edge generator occurring in a
square.

However, the second braid in \eqref{iop} does not have this form.
Thus it is the only reduced positive braid representation of
$16_{559253}$.

Now the diagrams of the closures of the braids in \eqref{iop} are
easily seen not to be transformable by mutations (e.g. switch them
to become alternating and calculate the Alexander polynomials~--
they differ). Also, since any other reduced positive braid diagram
of $16_{473363}$ must have $>16$ crossings, it cannot be a mutated
version of the 16 crossing positive braid diagram of $16_{559253}$,
and the two knots have no mutated positive braid diagrams.

There is one further braid positive mutant in this group, $16_{488722}$,
with 3 positive braid diagrams of 16 crossings, which are also not
mutated versions of the one of $16_{559253}$.

(Also the braid closures of the braids in \eqref{iop} are not mutants
in the complement of their braid axis, as the skein polynomials of
the 2-components links made up of braid closure and axis differ.)

Thus the problem to decide whether positive braids are mutants is
not solvable in the suggestive way.
\end{exam}

{\begin{figure}[htb]
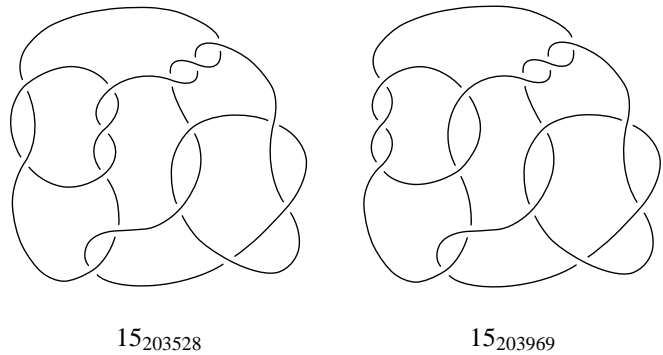

\[
\begin{array}{c@{\qquad}c}
\epsfsv{4cm}{t1-15_203528} & \epsfsv{4cm}{t1-15_203969}\\[19mm]
\ry{1.5em}15_{203528} & 15_{203969}
\end{array}
\]
\caption{These two mutants form with the knots on figure
\reference{figgc} a group of iterated mutants, but, contrarily
to the previous two knots, are braid positive. \label{figgd}}
\end{figure}
}

The next example, shows that~-- perhaps even worse~-- the property to
be braid positive it not preserved under mutation.

\begin{exam}
The knots $15_{203528}$ and $15_{203969}$ on figure \reference{figgd}
belong to the same group of (iterated) mutants as $15_{203432}$ and
$15_{203777}$  on figure \reference{figgc}, but their diagrams
are positive braid diagrams, so that they have positive braid
representations.
\end{exam}

\subsection{Computational details\label{cd}}

For the interested reader, here we make some supplementary remarks
on more of the details how the above examples were obtained
and checks were performed.

\subsubsection{\label{xx}Theorem \reference{thdn}}

First we address the proof of part \reference{item1}) of theorem
\reference{thdn}.
We wanted to check that if $\bt$ is a positive braid (word)
which is irreducible, i.e. \em{inter alia} modulo YB relations and
cyclic permutations not transformable into a word with an isolated
generator (i.e. letter appearing only once), then the subword $\bt_1$
of `$1$', `$2$' and `$3$' in $\bt$ has length at least 10. (We work
with words of integers 
`$i$' representing the corresponding generators $\sg_i$, and in the
sequel call these integers, despite being numbers, letters.) We already
know that the number of occurrences of `$1$' and `$3$' is at least $2$,
and of `$2$' is at least $4$, and that the word $\bt_1$ must have length
at least 9. To show is that no word of length 9 is possible.

This is done in 3 steps.
\begin{enumerate}
\def\labelenumi{\theenumi)}
\item All words of letters `$1$', `$2$' and `$3$' of length 9
(representing candidates for $\bt_1$) are generated. 
\item Irreducibility is tested. Irreducibility implies a number
of conditions on the word $\bt_1$. To reduce the number of cases,
we consider only words with maximal digit($=$ generator index) sum,
and among them only those, which are lexicographically minimal up
to cyclic permutations. For such words the number of occurrences of
`$1$', `$2$' and `$3$' is tested, and that the closure is not composite.
\item Connectedness test. Many braids $\bt$ can be discarded,
as the non-connectedness of their closure can already be seen from
their subwords $\bt_1$. We apply repeatedly to $\bt_1$ YB relations,
cyclic permutations,
and eliminate squares (pairs of consecutive copies) of `$1$' and `$2$'
(but not of `$3$', as there might be some letter `$4$' in $\bt$ between
these two copies we have discarded building $\bt_1$). Whenever this
procedure completely eliminates one of the letters `$1$' or `$2$',
the closure is not a knot, and the braid (word) can be discarded.
\end{enumerate}

These 3 checks already discard all possible
words of `$1$', `$2$' and `$3$' of length 9.

\subsubsection{\label{ST3}Proposition \reference{pdn} and theorem \reference{pdo}}

For examining 17 crossing 6-braids, it suffices again to consider
irreducible braids. To generate them, we use a similar
method. We split the word (of letters `$1$' to `$5$') into its
subwords $\bt_1$ of letters `$1$' to `$3$' and $\bt_2$ of letters `$3$'
to `$5$'. The words $\bt_1$ can be generated as above, and $\bt_2$ in a
similar way with the following modifications/remarks.

\begin{itemize}
\item Replace $1\to 5$ and $2\to 4$. (This replacement means that,
in composing $\bt_1$ and $\bt_2$ to $\bt$, we take $\bt$'s
maximal word representation modulo cyclic permutations, in which
the letters `$3$' are weighted higher than `$2$' and `$4$', and they
in turn are weighted higher than `$1$' and `$5$'.)
\item The connectedness check can be applied to $\bt_1$ and
$\bt_2$ in the same way as before (here to determine connectedness,
cyclic permutations are allowed also in $\bt_2$).
\end{itemize}

We can (up to flipping $2\lfra 4$ and $1\lfra 5$) assume that
$\bt_1$ is not shorter than $\bt_2$. Then for a 17 crossing irreducible
braid there are 3 possibilities, when keeping in mind the minimal number
of occurrences of the letters (at least $4$ for `$2$' and `$4$',
at least $2$ for `$1$', `$3$' and `$5$') and that $[\bt_i]\ge 10$.

\[
\begin{array}{c*{2}{|c}}
[\bt_1] & [\bt]_3 & [\bt_2] \\[0.3em]
\hline
\ry{1.2em} 10 & 3 & 10\\[0.3em]
10 & 4 & 11 \\[0.3em]
11 & 5 & 11
\end{array}
\]

Given the words $\bt_{1,2}$, there is a canonical way of
putting them together to obtain $\bt$: the (sub)words $w_{1,i}$ of `$1$'
and `$2$' between 2 occurrences of `$3$' in $\bt_{1}$ need to be
composed with the subwords $w_{2,i}$ of `$4$' and `$5$' between the
same occurrences of `$3$' in $\bt_{2}$. As the letters
`$1$' and `$2$' commute with `$4$' and `$5$', by
concatenating $w_{1,i}$ and $w_{2,i}$ we obtain the only relevant word.

The resulting diagrams can be checked for connectedness
(note that the above connectedness tests were just partial)
and the corresponding diagrams identified (they are about 1000).
It turned out, that only 6 knots occurred (two of 14 and four of 15
crossings), which all had braid positive minimal crossing diagrams.

To generate diagrams of $\le 16$ crossings, a different approach
was taken. From the alternating knot tables of \cite{KnotScape} the
fibred knots of the desired genus were selected by verifying the
(degree and leading coefficient of their)
Alexander polynomial. All flyped versions of their table diagrams
we generated (flypes are the moves of \cite{MenThis}),
and those knots where chosen which have
an alternating braid diagram (the property of a diagram to
be a braid diagram is not preserved by flypes, consider e.g.\ $7_7$).
Then the table diagrams of these knots were switched to be positive
(this commutes with flypes, so which alternating diagram of the knot is
taken is no longer relevant), and the resulting knots were identified.
These representations were found by examining the Jones polynomial.

A similar method to that of generating the `$3$' to `$5$' subwords was
used in the quest for (almost positive 17 crossing) 4-braid
representations of the knots in example \reference{exmbr}. The
word can be cyclically adjusted so as the negative generator to
be first. Then the output of the 3-5 subword program (with `$4$'
replaced by `$2$' and `$5$' replaced by `$1$'; relevant
is here that no cyclic permutations are
allowed to maximize the word) is appended to the negative crossing,
and connectedness (of the closure) checked for the result. (Here
the special meaning of the `$3$' in the connectedness check can be
eliminated.)

%
%
%
%
%
%
%
%
%
%
%
%
%
%
\def\hra{\hookrightarrow}
\def\sgn{\text{\operator@font sgn\,}}
\def\an#1{\left\langle#1\right\rangle}
\def\BR#1{\left\lceil#1\right\rceil}

\section{Braid index inequalities}

\subsection{On a 4-braid criterion of Jones}

In the famous paper \cite{Jones}, where Jones announced his fundamental
discovery of a relation between $C^*$-algebras and Markov traces on
braid groups, he gave also some results concerning applications of his
new invariant to braids. Most of these, and many more, results have
subsequently appeared with proof in his work-out \cite{Jones2}.
One of these results (theorem 22 in \cite{Jones}) was a formula
relating the Jones $V$ and Alexander polynomial $\Dl$ (normalized
so that $\Dl(1)=1$ and $\Dl(t^{-1})=\Dl(t)$) of 4-braid knots.

\begin{theorem}(\cite[proposition 11.11]{Jones2})
If a knot $K$ is the closure of a 4-braid $\bt$ of exponent sum $[\bt]$,
then
\begin{eqn}\label{Je}
t^{-[\bt]}V_K(t)+t^{[\bt]}V_K(1/t)\,=\,\bigl(t^{-3/2}+t^{-1/2}+t^{1/2}+t^{3/2}
\bigr)(t^{-{[\bt]}/2}+t^{{[\bt]}/2})-\bigl(t^{-2}+t^{-1}+2+t+t^2\bigr)\Dl_K(t)\,.
\end{eqn}
\end{theorem}

As a consequence of this theorem, in \cite{Jones} Jones announced an
obstruction to braid index 4 for knots (corollary 24), namely that the
value $\Dl(e^{2\pi i/5})$ for such knots must be of norm at most 6.5.
This
result did not appear with proof in \cite{Jones2}, and my attempts to
recover it failed. Finally, it turned out that the result, as stated
there, is in fact not correct, and there is a counterexample.

\begin{exam}
The knot $13_{9221}$ of \cite{KnotScape} shown on figure \reference{%
fig2br} has a (braid) diagram with 4 Seifert circles, and thus
braid index at most $4$. (We have $\spn_v P=6$, where $P$ is the HOMFLY
polynomial \cite{HOMFLY}, 
so by MWF \cite{WilFr,Morton} the braid index is indeed $4$). However,
its Alexander
polynomial is $\Dl(t)=t^{-3}-10t^{-2}+29t^{-1}-39+29t-10t^2+t^3$, which
at $t=e^{2\pi i/5}$ evaluates to $19\sqrt{5}-49\approx -6.5147084$.
\end{exam}

{
\begin{figure}[htb]
\[
\begin{array}{c}
\epsfsv{4cm}{t1-13-9221} \\[2.2cm]
13_{9221}
\end{array}
\]
\caption{\label{fig2br}}
\end{figure}
}

The value is still very close to $6.5$, and as this constant does
not appear very natural, it is suggestive that it might have been
obtained by rounding (possibly erronously the \em{difference} of the
estimates of \eqref{2} and \eqref{2.5} given below was taken). However,
there is some evidence that the bound
cannot be fixed even just by a minor improvement. Instead we present a
criterion with a larger, but definitely correct bound.

\begin{proposition}\label{pp}
If a knot $K$ has $\big|\Dl_K(e^{2\pi i/5})\big|>6+2\sqrt{5}
\approx 10.472136$, then $K$ is not a closed $4$-braid.
\end{proposition}

\proof The formula \eqref{Je} for $t=e^{2\pi i/5}$ simplifies
to
\[
t^{-{[\bt]}}V_K(t)+t^{[\bt]}V_K(1/t)\,=\,2\cos {{[\bt]}\pi /5}-\Dl_K(t)\,,
\]
which, as $t$ has unit norm, and thus $V_K(1/t)=\overline{V_K(t)}$,
gives
\begin{eqn}\label{2}
\Dl_K(t)\,=\,2\cos {{[\bt]}\pi /5}-2\Re\,(t^{-{[\bt]}}V_K(t))\,.
\end{eqn}
We have 
\begin{eqn}\label{2.5}
\big|2\cos {{[\bt]}\pi /5}\big|\le 2
\end{eqn}
and from \cite[proposition 14.6]{Jones2} also
\begin{eqn}\label{3}
2\Re\,(t^{-{[\bt]}}V_K(t))\,\le\,2\big|V_K(e^{2\pi i/5})\big|\le 16\cos^3
(\pi /5)\,=\,\frac{(1+\sqrt{5})^3}{4}\,=\,4+2\sqrt{5}.
\end{eqn}
Putting \eqref{2.5} and \eqref{3} into \eqref{2} gives the result. \qed

By remarking that $10\nmid [\bt]$ (as $[\bt]$ is odd for a knot),
one can slightly improve the upper estimate in \eqref{3}, obtaining that
\begin{eqn}\label{eq3}
\Dl_K(e^{2\pi i/5})\le\frac{9+5\sqrt{5}}{2}\approx 10.09017
\end{eqn}
for a $4$-braid knot $K$.

The inequality \eqref{2} is clearly sharp, for $5\mid {[\bt]}$ (it is
not surprising that in the above example indeed ${[\bt]}=5$). Also,
the second estimate in \eqref{3} is trivially sharp for \em{links} of
braid index at most $4$ (take the 4 component unlink), but
the denseness result in \cite[proposition 14.6]{Jones2} was sharpened in
\cite{gen2} to show that it remains true even if one restricts his
attention to knots, so that the only way to improve the bound
in proposition \reference{pp} (resp. \eqref{eq3}) along these lines of
argument is to improve the left inequality in \eqref{3}, that is, to
show that $t^{-{[\bt]}}V_K(t)$ is sufficiently far from the real line.
This appears, however, unlikely as well, and thus a much better
constant than the one given above can probably not be obtained.

The reason why Jones's 4-braid criterion never attracted particular
attention is possibly that briefly later the (much more effective)
MWF inequality was found. Indeed, for non-alternating prime knots
of at most 14 crossings $13_{9221}$ was the only one which violated
Jones's original (and insufficient) condition on $\Dl(e^{2\pi i/5})$,
and which had $\spn_v P\le 6$. An advantage of the (corrected)
criterion involving $\Dl$ remains, however, that it is applicable also
to very complicated knots because of the polynomial complexity of $\Dl$.

Finally, we mention that the above knot, $13_{9221}$, has another
interesting property which is discussed in a joint paper with Mark
Kidwell \cite{KS}.

\subsection{On the 2-cabled MWF inequality\label{S3}}

In \cite{MorSho1}, Morton and Short introduced a way to circumvent
MWF's failure to estimate sharply the braid index of a knot $K$
by applying the inequality on a 2-cable $K_2$ of $K$, a satellite
around $K$ with a pattern intersecting each meridian disc of the
solid torus twice and in the same direction. As this is a 2-braid,
the satellite $K_2=K_{2,w}$ is uniquely determined by the writhe $w$
of this braid (the satellite is connected or disconnected depending
on the parity of $w$). To obtain a braid representation of $K_{2,w}$
from a braid representation $\bt$ of $K$,
a generator $\sg_i$ in $\bt$ is replaced by
$\sg_{2i}\sg_{2i-1}\sg_{2i+1}\sg_{2i}$, and the result is multiplied
by $\sg_1^{-2[\bt]+w}$. Thus $b(K_{2,w})\le 2b(K)$, and applying MWF
on $K_{2,w}$ we obtain (for any $w\in\bZ$)
\[
b(K)\ge\BR{\frac{\spn_vP_{K_{2,w}}/2+1}2}.
\]
This inequality completely determines the braid index of the Rolfsen
knots, for which MWF failed itself.
However, even this inequality sometimes fails, as shows the
example below. As for this example the braid index to exclude is $4$
(and hence the 3-braid conditions of \cite{Murasugi3} and \cite{Murakami}
are not relevant), this gives an example of a knot, on which all
previously applied braid index criteria fail to estimate sharply
the braid index.

\begin{exam}
The knot $14_{45759}$ is depicted on figure \reference{fig3br}.
Its $P$ polynomial is shown in table \reference{PP}, and estimates the
braid index to be at least 3. This is however seen not to be exact
from the $P$ polynomial of a 2-cable knot of $P$ shown below in
table \reference{PP}, estimating the braid index to be at least 4.
However, even this estimate is not exact. To see this, we use that
the knot is achiral and need to go a little behind the MWF inequality.
This inequality was the consequence of the following two inequalities.

\begin{theorem}(\cite{Morton,WilFr})\label{TWF}
If $K=\hat\bt$, $\bt\in B_n$, then $\md_vP_K\ge [\bt]-n+1$ and
$\Md_vP_K\le [\bt]+n-1$.
\end{theorem}

If $14_{45759}$ were a closed $4$-braid $\bt$, then from the
$P$ polynomial we see that only $[\bt]=\pm 1$ can occur, and
indeed both values do because of achirality. Call these
braids $\bt_{1,2}$. But then, taking the 2-cable of $\bt_{1,2}$
we obtain $8$-braid representations of $K_2$ of different writhe.
However, the fact that MWF is sharp for these $8$-braids, contradicts
one for the 2 inequalities above. Thus the braid index of $14_{45759}$
is at least $5$. 
\end{exam}

{
\begin{figure}[htb]
\[
\begin{array}{c}
\epsfsv{4cm}{t1-14-45759} \\[2.3cm]
14_{45759}
\end{array}
\]
\caption{\label{fig3br}}
\end{figure}
}

\begin{table}
{\small
\begin{verbatim}
14 45759     0     8
         -2   2          -3    -5    -3
         -2   2           1     5     1
         -2   2           2     3     2
         -2   2          -1    -4    -1
          0   0                 1


57 45759     0    24
         -4   6                    28      84      78      -5     -45     -19
         -6   8           -35    -236    -527    -394     120     350     183      33
         -6   8           405    1342    1658     665    -578   -1157    -892    -276
         -6   8         -1423   -4237   -3467    -167    1181    2092    2276     759
         -6   8          2303    7563    5679    -166   -1245   -2672   -3300    -937
         -6   8         -1965   -8202   -7643   -1601     706    2930    2892     583
         -6   8           936    5755    7904    4143    -215   -2610   -1625    -190
         -6   8          -250   -2713   -5823   -4401      33    1624     597      31
         -6   8            35     857    2902    2583      -2    -640    -138      -2
         -6   6            -2    -173    -940    -901       0     150      18
         -4   6                    20     188     186       0     -19      -1
         -4   4                    -1     -21     -21       0       1
         -2   0                             1       1
\end{verbatim}
}
\caption{The $P$ polynomial of $14_{45759}$ and a 2-cable knot
of it.\label{PP}}
\end{table}

\subsection{The Jones conjecture\label{S7}}

The example and the reasoning applied in the previous section
can be possibly made more general.

A $n$-parallel $K_\gm$ of $K$ is a satellite around $K$ of zero
framing with pattern being a closed $n$-string braid $\gm$ in the
solid torus given by the complement of its braid axis.	

\begin{lemma}
If $\bt\in B_n$ is a braid representation for $K$, then $\an{\bt}_\gm
\in B_{kn}$ is a braid representation for $K_\gm$, $\gm\in B_k$.
Here 
\begin{eqn}\label{***}
\an{\bt}_\gm\,=\,\gm\cdot\{\bt\}^k\cdot \Dl_k^{-2[\bt]}\,,
\end{eqn}
where $\{\bt\}^k$ is obtained from $\bt$ by the
replacement 
\[
\sg_i^\eps\,\longmapsto\,\prod_{n=1}^{2k-1}\,\prod_
{j=(ik-\min(n,2k-n)+1)/2}^{(ik+\min(n,2k-n)-1)/2}\,\sg_{2j}^\eps\,,
\qquad\eps=\pm 1\,,
\]
and
\[
\Dl_k\,=\,\prod_{j=1}^{k-1}\,\prod_{l=1}^{k-j}\,\sg_l
\]
is the square root of the center generator of $B_k$
(the first and third factor on the right of \eqref{***}
are meant w.r.t. the inclusion $B_k\hra B_{nk}$).
\end{lemma}

\proof This is a well-known and trivial fact (although seldom
stated in such explicity). Taking the diagram $\hat\bt$
of $K$, we add $\big|\,[\bt]\,\big|$ kinks of sign $-\sgn[\bt]$,
cable the diagram (under which $\bt$ is taken to $\{\bt\}^k$),
and remove the kinks, obtaining $\Dl_k^{-2\sgn[\bt]}$ for each
kink. \qed

Considering a (connected or disconnected) parallel $n$-cable knot
(or link) $K_n$ of a knot $K$ (the choice of $\gm\in B_n$ is
no longer relevant), we can apply MWF to $K_n$ and
use $b(K_n)\le nb(K)$, thus obtaining an infinite series
of inequalities (for any $n$)
\begin{eqn}\label{st}
b(K)\ge\BR{(\spn_vP_{K_n}/2+1)/n}.
\end{eqn}
The practical problem with these inequalities is that the calculation
of $P_{K_n}$ is impossible for $n\ge 3$ and any, even moderately
interesting, example $K$. Nevertheless, one can ask whether \eqref{st}
can be made unsharp for more small values of $n$, or even for all
$n$. 

It turns out that this problem is related to
one of the still unsolved conjectures made by Jones
briefly after his discoveries. It is as follows.

\begin{conj}(Jones \cite{Jones2})\label{cj1}
If $\bt,\bt'\in B_{b(L)}$ and $\hat\bt=\hat\bt'=L$, then $[\bt]=[\bt']$.
\end{conj}

We include a brief historical review. The conjecture was first very
implicitly mentioned in Jones's paper \cite{Jones2}. Later,
some main publicity to it was given by Birman in her paper
with Menasco \cite{BirMen}, where it was proved (corollary p.\ 267) that
each link has at most finitely many writhes of minimal braid 
representations. However, this was proved previously by Morton in his
paper \cite{Morton}, and also in \cite{WilFr}, in a much less
sophisticated way, not only for minimal, but for any arbitrary
fixed strand number, and in stronger form, with very explicit
lower and upper bounds to the writhe in terms of the degrees of the
skein polynomial. In particular the Jones conjecture follows to be true
for links with sharp MWF inequality, or in fact for any link for which
the braid index can be determined by applying MWF on some parallel
cable, as in \cite{MorSho2}. Since this relationship
was apparently not previously realized, it will be explained below.

Also, another statement of Birman needs correction.
It is claimed that the Jones conjecture ``is
known to be true'' for 4-braids ``by the work of Jones'' \cite{Jones2}.
However, at least in the case that a braid index 4 knot, as $10_{132}$,
has the same skein polynomial as some knot of smaller braid index,
Jones's formulas (\cite[\S 8]{Jones2}) will certainly
not be able to exclude multiple writhes of 4-braid representations,
and already to classify when such duplications of the skein
polynomial occur seems impossible. At least this was never (to
the best of my knowledge) carried out by Jones. Jones himself only
wrote that his formulas ``lend some weight ot the possibility''
(\cite[p. 357 bottom]{Jones2}) this conjecture to be true.
Thus the conjecture should be considered open in
the 4-strand case (for 3-braids it was settled in \cite{BirMen2}).
In fact, a certain importance of 4-braids for this conjecture
will be established later.

An interesting special case of the Jones conjecture was
addressed in a question raised independently by P.\ Johnson:

\begin{question}
Is there an achiral knot $K$ of even braid index?
\end{question}

The argument used in the previous example immediately
shows that such a knot would be a counterexample to
the Jones conjecture.

A generalization of our argument in \S\reference{S3} shows
the following explicit version of the Birman-Menasco result:

\begin{theorem}\label{>>>}
For any knot $K$ and any $k\ge b(K)$ we have
\begin{eqn}\label{opq}
k-b(K)+1\le d_{k,K}=\#\{\,[\bt]\,:\,\bt\in B_k,\ \hat\bt=K\,\}\le
k-\max_{n,K_n}\left(\mbox{r.h.s. of \eqref{st}}\right)\,+1\,.
\end{eqn}
\end{theorem}

Note, that conjecture \reference{cj1} is equivalent to
$d_{k,K}=k-b(K)+1$ for any $K$ and $k\ge b(K)$.

\proof[of theorem \reference{>>>}] First observe that
if $\bt,\bt'\in B_n$, and $\gm\in B_k$, then
\begin{eqn}\label{lml}
\big|\,[\an{\bt}_\gm]-[\an{\bt'}_\gm]\,\big|\,=\,k
\big|\,[\bt]-[\bt']\,\big|\,.
\end{eqn}

To see this, use $[\Dl_k]=\frac{k(k-1)}{2}$ and $[\{\bt\}_k]=k^2[\bt]$.
The result is straightforward from \eqref{***}. 

Now, it is immediate from the inequalities of theorem \reference{TWF}
that if $\bt,\bt'\in B_n$ with $\hat\bt=\hat\bt'=K$, then
\begin{eqn}\label{lmm}
\frac{1}{2}\spn_vP_K+1\le n-\frac{\big|\,[\bt]-[\bt']\,\big|}{2}\,.
\end{eqn}


This already shows that the sharpness of the MWF inequality implies
the truth of the Jones conjecture. However, we can get this relationship
now in a more general version.

Applying \eqref{lmm} on $\an{\bt}_k$ and
$\an{\bt'}_k$, and using \eqref{lml}, we see that
if $\bt,\bt'\in B_n$ with $\hat\bt=\hat\bt'=K$, and $\gm\in B_k$, then
\begin{eqn}\label{lmn}
\frac{1}{k}\left[
\frac{1}{2}\spn_vP_{K_\gm}+1\right]
\le n-\frac{\big|\,[\bt]-[\bt']\,\big|}{2}\,.
\end{eqn}


Let $D_{k,K}:=\{\,[\bt]\,:\,\bt\in B_k,\ \hat\bt=K\,\}$. As 
\[
d_{k,K}=\#D_{k,K}\le
\max_{[\bt],[\bt']\in D_{k,K}}\frac{\big|\,[\bt]-[\bt']\,\big|}{2}+1\,,
\]
we get from \eqref{lmn} the second inequality in
\eqref{opq}. The first inequality is trivial (take a minimal
braid representation and stabilize in all possible ways). \qed

The important consequence is the case $k=b(K)$ and $d_{k,K}\ge 2$:

\begin{corr}
If $K$ is a counterexample to conjecture \reference{cj1}, then
\eqref{st} is unsharp for $K$ for any choice of $n$ and $K_n$. \qed
\end{corr}

This means that, provided we want to give a counterexample
to conjecture \reference{cj1} and even have found $\bt$ and
$\bt'$, we cannot prove their minimality using any
of the inequalities \eqref{st}. This shows why the quest for
alternatives to MWF is worthwhile. There are such criteria,
due to Jones \cite{Jones2}, but (in particular, because MWF
performs well very often) it is difficult to find examples
where these criteria show more powerful than MWF. Although
such examples exist (two are shown in the next section), they are
very rare. Another handy criterion
is Murakami's 3-braid formula \cite[corollary 10.5]{Murakami}. This
criterion is very effective~-- for example it excludes (without a
single failure!) from having
braid index 3 all 916 non-alternating prime 15 crossing knots
with MWF bound 3. (The candidates for exponent sums of 3-braids being
taken from the skein polynomial.) As it applies only for 3-braid knots,
at present it seems unrealistic to find a counterexample to the Jones
conjecture in braid index $>4$. For the more optimistic
readers, the corollary can also be taken as evidence for the
conjecture.

In this situation I initiated a large
computer experiment. I selected knots with MWF bound $\le 3$ from the
tables of \cite{KnotScape}, for which the Murakami test excludes
braid index 3, but for which \eqref{3} is satisfied for the (or at
least) two values of $[\bt]$ for $\bt\in B_4$ admitted by $P$
from the inequalities in theorem \reference{TWF}. Then I calculated $P$
of a 2-cable of these knots. This was already a non-trivial task.
For knots for which I could find minimal crossing number diagrams
with $\le 5$ Seifert circles, I applied Vogel's algorithm \cite{Vogel}.
The simplest (=lowest crossing number) braid representation 
obtained was 2-cabled (a generator $\sg_i$ replaced by
$\sg_{2i}\sg_{2i-1}\sg_{2i+1}\sg_{2i}$) and processed by the program of
\cite{MorSho1,MorSho2}. For the other knots the DT notation \cite{%
DT} of a 2-cable knot was generated from the DT notation of the knot
(given in the tables of KnotScape), and the polynomial calculation
program of KnotScape (a variation of the Millett-Ewing
program) was used. This way (and with some assistance of Ken Millett
for the hardest examples) I verified all prime knots up to 15 crossings,
and many of the knots of 16
crossings, to which one of the programs of Millett-Ewing or
Morton-Short program was applicable 
(for some knots both programs failed
due to memory and time constraints). Although it did
not give a counterexample to the Jones conjecture, this experiment
found examples like the one of \S\reference{S3}.

The following corollary summarizes for which classes of knots
the sharpness of MWF, and hence the truth of the Jones
conjecture, is known (see \cite{WilFr,Murasugi2}).

\begin{corr}
The Jones conjecture is true for
\begin{itemize}
\item alternating fibred knots,
\item rational knots, and
\item positive braid knots with a full twist (i.e., closures of
braids of the form $\Dl_n\ap\in B_n$ with $\ap$ positive).
\end{itemize}
\end{corr}

\subsection{The Jones polynomial at roots of unity\label{SU}}

One of Jones's original criteria for the braid index
came from the positivity of (a scalar product on) a $C^*$-algebra
related to the values of his polynomial at primitive roots of unity.

\begin{theorem}(\cite[proposition 14.6]{Jones2})
If a knot $K$ has an $n$-braid representation and $k\ge 3$, then
\begin{eqn}\label{jc}
\big|\,V_K\bigl(e^{2\pi i/k}\bigr)\,\big|\,\le\,
\bigl(2\cos \pi/k\bigr)^{n-1}\,.
\end{eqn}
\end{theorem}

{
\begin{figure}[htb]
\[
\begin{array}{c@{\qquad}c}
\epsfsv{4cm}{t1-13-9465} & \epsfsv{4.2cm}{t1-14-41800} \\[2.2cm]
13_{9465} & 14_{41800}
\end{array}
\]
\caption{\label{figur}}
\end{figure}
}

\begin{table}[htb]
{\small
\begin{verbatim}
13 9465     0    10
         10  14         -11   -15    -5
         10  14          40    35     5
         10  14         -57   -28    -1
         10  12          36     9
         10  12         -10    -1
         10  10           1

13    9465        5   15         1    0    1    0    1   -1    0   -1    0   -1    1
\end{verbatim}
}
\caption{The $P$ and $V$ polynomial of $13_{9465}$.\label{PV}}
\end{table}

This criterion turned out to be of less practical relevance
than MWF, which is much more direct to apply and often more efficient.
That, however, the inequalities \eqref{jc} can sometimes give better
estimates, and thus need to be rehabilitated, shows the following
example.

\begin{exam}
The trefoil cable knot $K=13_{9465}=(3_1)_{2,7}$, shown on figure
\reference{figur}, has a $P$ polynomial giving the MWF bound $3$
(this was noticed already in \cite{WilFr}). However, when considering
its Jones polynomial (see table \reference{PV}) and using \eqref{jc}
for $k=10$, we find
\[
V_K\bigl(e^{\pi i/5}\bigr)\,=\,-2-\frac{3+\sqrt{5}}{2}\sqrt{
\frac{5-\sqrt{5}}{2}}i\,,
\]
and thus
\[
\big|\,V_K\bigl(e^{\pi i/5}\bigr)\,\big|^2\,=\,9+2\sqrt{5}\approx
13.472136\dots\,,
\]
which exceeds
\[
\bigl(2\cos \pi/10\bigr)^4=\frac{15+5\sqrt{5}}{2}\approx 13.0901699\dots
\]
Thus $13_{9465}$ cannot be a 3-braid, and has braid index $4$
(the diagram shows it to be a 4-braid). Similar examples
(which are not closures of positive braids) are $14_{41800}$,
$15_{138678}$, $15_{141289}$ and $15_{251434}$.
\end{exam}

The quest for such examples was motivated by the question
(raised by Birman in problem 10.1 of \cite{Morton4}) on the
realizability of the MWF bound among knots of given $P$ polynomial.
This question is thus answered to the negative, first for
MWF bound $3$, but examples for all higher values of the bound can be
constructed by iteratedly taking (the polynomials of)
connected sums of the above knot $13_{9465}$
with itself or with $5_1$ (or their mirror images). The cases
remaining open are for MWF bound $1$ and $2$.

\begin{question}
Are there $P$ polynomials of $v$-span $\le 2$ other than those
of the $(2,n)$-torus knots and links?
\end{question}

The problem for specific examples of polynomials is still difficult.
Among the polynomials of the $5$ Rolfsen knots with unsharp MWF
(the two 9 crossing knots are addressed in problem 10.2 of
\cite{Morton4}), the MWF bound is realizable for the polynomials of
$10_{156}$ and $10_{132}$ (by $8_{16}$ and $5_1$ respectively,
see \cite[notes p.\ 386]{Jones2}), but the status of the
remaining 3 polynomials is undecided. There is no prime knot
of $\le 16$ crossings duplicating the $P$ polynomial of $9_{42}$,
and only one, $16_{730458}$, duplicating this of $9_{49}$. For the
polynomial of $10_{150}$ there are two duplications in the tables,
the knots being $13_{4977}$ and $13_{6718}$. However, for all these
3 duplicating knots the Murakami formula showed that they are
not of braid index 3 (for 2 of them the 2-cabled MWF inequality
even showed that the braid index is at least 5).

\section{The genus and HOMFLY polynomial}

A final, and somewhat unrelated, collection of examples
concerns two other conjectured relations of the degrees of the
skein polynomial, this time related to genera.

\subsection{Morton's conjecture}

Briefly after the discovery of the skein polynomial,
Morton \cite{Morton4} posed the question whether for
any link $L$, 
\begin{eqn}\label{mc}
1-\chi(L)\,\ge\,\md_vP(L)\,,
\end{eqn}
where $\chi(L)$ is the maximal
Euler characteristic of a spanning surface for $L$ (if $L$ is
a knot, then $1-\chi(L)$ is twice the genus of $L$).

The motivation for this question was the fact that both hand-sides
of \eqref{mc} are estimated below by $[\bt]-n(\bt)+1$ for any
braid $\bt$ with $\hat\bt=L$. For the r.h.s. this was, as noted,
proved by Morton himself, while for the l.h.s. it is a famous
inequality of Bennequin \cite[theorem 3]{Bennequin}.

Morton's conjecture resisted solution for a long time. \eqref{mc}
has been proved for homogeneous links \cite{Cromwell}, 3-braid
links \cite{DasMan}, and verified computationally for prime
knots up to 16 crossings \cite{deg}. (Since both hand-sides of
\eqref{mc} are additive under connected sum, it suffices to check
prime links.) Also, further attempts have been made \cite{Cromwell3}.

Here we settle Morton's conjecture negatively by means of
counterexamples.

{\begin{figure}[htb]
\[
\begin{array}{c}
\epsfsv{5.7cm}{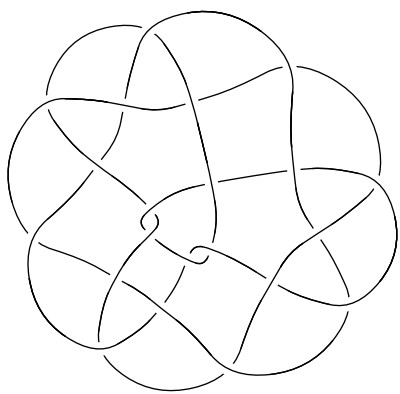}
\end{array}
\]
\caption{A counterexample to Morton's conjecture.\label{figmc}}
\end{figure}
}

\begin{exam}
Consider the knot $K$ on figure \reference{figmc}. 
It has the following 4-braid representation:
\[
\sg_1^2 (\sg_1 \sg_2 \sg^{-1}_1) (\sg_2 \sg_3 \sg^{-1}_2) (\sg_1
\sg_2 \sg^{-1}_1) (\sg_2 \sg_3 \sg^{-1}_2) \sg_3 \sg_1 (\sg_1
\sg_2 \sg^{-1}_1) (\sg_1 \sg_2 \sg_3 \sg^{-1}_2 \sg^{-1}_1) \sg_3
\]
From this representation it is evident that $K$ has a genus 4
surface obtained by connecting the 4 discs of the strands
by the 11 bands indicated by the parenthesized subwords (see \cite{%
Rudolph2}). That this Seifert surface has minimal genus follows from
Bennequin's inequality, since all the bands are `positive'.

Thus $K$ has genus 4. However, a calculation shows that $\md_vP(K)=10$.
\end{exam}

There are 7 further such examples, given by the 4-braids:
\begin{verbatim}
1 1 1 2 -1 2 1 3 1 2 -1 2 2 3 -2 1 2 -1 2 3 -2 
1 1 1 2 -1 2 1 3 1 2 -1 2 3 -2 1 2 -1 1 2 3 -2 -1 2 3 -2 
1 1 1 2 -1 2 2 3 -2 1 2 -1 2 3 -2 3 1 2 -1 1 2 3 -2 3 -1 
1 1 1 2 -1 2 1 3 2 -1 1 2 3 -2 -1 2 3 -2 1 1 2 -1 2 3 -2 
1 1 1 2 -1 2 1 3 2 -1 1 2 3 -2 -1 2 3 -2 1 2 -1 2 3 -2 3 
1 1 1 2 -1 2 3 -2 1 1 2 -1 1 2 3 -2 -1 2 3 -2 1 2 -1 1 2 3 -2 3 -1 
1 1 1 2 -1 2 3 -2 1 2 -1 1 2 3 -2 -1 2 3 -2 3 1 2 -1 1 2 3 -2 3 -1 
\end{verbatim}

For all 8 examples, KnotScape manges to reduce the diagrams
to 21 crossings, while from the degree $\Md_zF=17$ of their
Kauffman polynomial one concludes that their crossing number
is at least 19. (They can be distinguished by the Kauffman
polynomial, although the skein polynomials of several of them
coincide.)

\subsection{Morton's canonical genus inequality}

\begin{exam}
A final example concerns another inequality
of Morton proved in \cite{Morton}: $\Md_zP_K\le \tl g(K)$.
Here $\tl g(K)$ is the canonical genus of a knot $K$, the minimal
genus of the canonical Seifert surfaces of all its diagrams
(see e.g. \cite{gen1,Kobayashi}). Another obvious inequality is
$g(K)\le \tl g(K)$. In a comparison of the 2 estimates for $\tl g$, in
\cite{Morton} it was remarked that knots exist with $2g<\Md_zP_K$.
Since Morton's inequality is exact for very many knots (in particular
all knots up to 12 crossings), in \cite{Morton5} I asked whether for
some knot the opposite relation $2g>\Md_zP_K$ can occur.
Such examples indeed exist, and were found by implementing
Gabai's method of disc decomposition \cite{Gabai} on
canonical Seifert surfaces of special diagrams. See figure
\reference{figMg}. This gives another example of knots with
unsharp Morton inequality, after the ones found in \cite{gen2}.
\end{exam}

{\begin{figure}[htb]
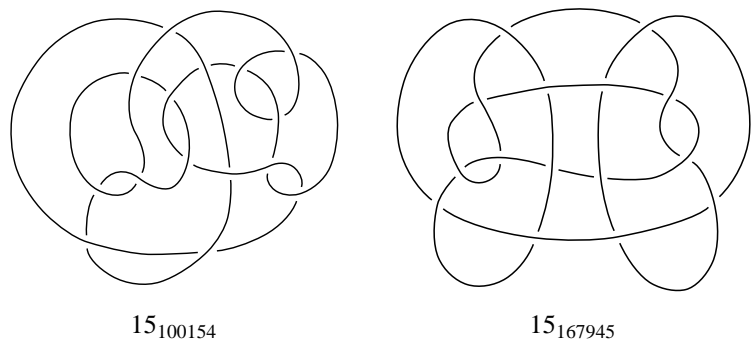

\[
\begin{array}{c@{\qquad}c}
\epsfsv{4cm}{t1-15_100154} & \epsfsv{4cm}{t1-15_167945}\\[3mm]
\ry{1.5em}15_{100154} & 15_{167945}
\end{array}
\]
\caption{Two knots for which $2g>\Md_zP$. In both cases
$\Md_zP=6$, while the canonical surfaces of the above
diagrams are disc decomposable, and hence $g=4$. There are
10 further 15 crossing examples of this type.\label{figMg}}
\end{figure}
}

\noindent{\bf Acknowledgement.} I would wish to thank to Ken Millett
for his help in doing a part of the calculations, and also to V.~Jones,
H.~Morton and K.~Murasugi for helpful remarks and discussions.

{\small
\let\old@bibitem\bibitem
\def\bibitem[#1]{\old@bibitem}

}

\begin{thebibliography}{48}
\bibitem[Ad]{Adams} C.~C.~Adams, \em{Das Knotenbuch}, Spektrum
        Akademischer Verlag, Berlin, 1995 (\em{The knot book},
        W.~H.~Freeman \& Co., New York, 1994).
\bibitem[Ad2]{Adams2} \bysame\ et al., \em{Almost alternating links},
        Topol. Appl. {\bf 46} (1992), 151--165.
\bibitem[Al]{Alexander} J. W.~Alexander, \em{Topological invariants
        of knots and links}, Trans.\ Amer.\ Math.\ Soc. {\bf 30} (1928),
        275--306.
\bibitem[An]{Aneziris} C.\ Aneziris, \em{The mystery of knots}. Computer
	programming for knot tabulation. Series on Knots and Everything
	{\bf 20}, World Scientific, 1999.
\bibitem[Be]{Bennequin} D.~Bennequin, \em{Entrelacements et \'equations
	de Pfaff}, Soc.\ Math.\ de France, Ast\'erisque {\bf 107-108}
	(1983), 87--161.
\bibitem[Bi]{Birman} J.~S.~Birman, \em{On the Jones polynomial of
	closed 3 braids}, Invent. Math. {\bf 81} (1985), 287--294.
\bibitem[BM]{BirMen} \bysame\ and W.~W.~Menasco, \em{Studying knots
	via braids VI: A non-finiteness theorem},
	Pacific J.~Math. {\bf 156} (1992), 265--285.
\bibitem[BM]{BirMen2} \bysame\ and \bysame, \em{Studying knots
	via braids III: Classifying knots which are closed 3 braids},
	Pacific J.~Math. {\bf 161}(1993), 25--113.
\bibitem[BW]{BirWil} \bysame\ and R.~F.~Williams, \em{Knotted
	periodic orbits in dynamical systems~- I, Lorenz's equations},
	Topology {\bf 22(1)} (1983), 47--82.
\bibitem[BLM]{BLM} R.~D.~Brandt, W. B. R. Lickorish and K. Millett,
	\em{A polynomial invariant for unoriented knots and links},
	Inv. Math. {\bf 74} (1986), 563--573.
\bibitem[Bu]{Busk} J. v. Buskirk, \em{Positive links have positive
	Conway polynomial}, Springer Lecture Notes in Math. {\bf 1144}
	(1983), 146--159.
\bibitem[CSV]{CSV} L. Carlitz, R. Scoville and T. Vaughan, 
	\em{Some arithmetic functions related to Fibonacci numbers},
	Fib. Quart. {\bf 11} (1973), 337--386.
\bibitem[Co]{Conway} J.~H.~Conway, \em{On enumeration of knots and
	links}, in ``Computational Problems in abstract algebra''
        (J.~Leech, ed.), 329-358, Pergamon Press, 1969.
\bibitem[Cr]{Cromwell} P. R. Cromwell, {\em Homogeneous links},
	J. London Math. Soc. (series 2) {\bf 39} (1989), 535--552.
\bibitem[Cr2]{Cromwell2} \bysame, {\em Positive braids are visually
	prime}, Proc. London Math. Soc. {\bf 67} (1993), 384--424.
\bibitem[Cr3]{Cromwell3} \bysame, {\em A note on Morton's conjecture
	concerning the lowest degree of a $2$-variable knot polynomial},
	Pacific J. Math. {\bf 160(2)} (1993), 201--205. 
\bibitem[CM]{MorCro} \bysame\, and H.~R.~Morton, \em{Positivity of knot
	polynomials on positive links}, J. Knot Theory Ramif. {\bf 1}
	(1992), 203--206.
\bibitem[DM]{DasMan} O. T. Dasbach and B. Mangum, \em{On McMullen's and
	other inequalities for the Thurston norm of link complements},
	Algebraic and Geometric Topology, {\bf 1} (2001), 321--347.
\bibitem[DT]{DT} C.~H.~Dowker and M.~B.~Thistlethwaite,
	\em{Classification of knot projections}, Topol. Appl.
	{\bf 16} (1983), 19--31.
\bibitem[Fi]{Fiedler} Th.~Fiedler, {\em A small state sum for knots},
	Topology {\bf 32~(2)} (1993), 281--294.
\bibitem[FW]{WilFr} J.\ Franks and R.\ F.\ Williams, \em{Braids and the
   	Jones-Conway polynomial}, Trans.\ Amer.\ Math.\ Soc. {\bf 303}
      	(1987), 97--108.
\bibitem[Ga]{Gabai} D.~Gabai, \em{Foliations and genera of links},
	Topology {\bf 23} (1984), 381--394.
\bibitem[H]{Hirzebruch} F.\ Hirzebruch, \em{Singularities and exotic
	spheres}, Seminaire Bourbaki {\bf 10}, Exp. No. 314, 13--32,
	Soc. Math. France, Paris, 1995.
\bibitem[H]{HOMFLY} P. Freyd, J. Hoste, W. B. R. Lickorish,
	K. Millett, A. Ocneanu and D. Yetter, {\it A new polynomial
	invariant of knots and links}, Bull. Amer. Math. Soc.
	{\bf 12} (1985), 239--246.
\bibitem[Ho]{Ho} C.~F.~Ho, \em{A polynomial invariant for knots and
	links~-- preliminary report}, Abstracts Amer. Math. Soc. {\bf
	6} (1985), 300.
\bibitem[HT]{KnotScape} J.~Hoste and M.~Thistlethwaite, {\em
        KnotScape}, a knot polynomial calculation and table access
	program, available at {\tt http://www.math.utk.edu/\~morwen}.
\bibitem[J]{Jones} V.~F.~R.~Jones, {\em A polynomial
        invariant of knots and links via von Neumann algebras},
        Bull. Amer. Math. Soc. {\bf 12} (1985), 103--111.
\bibitem[J2]{Jones2} \bysame, {\em Hecke algebra representations of
  	of braid groups and link polynomials}, Ann. of Math.
	{\bf 126} (1987), 335--388.
\bibitem[K]{Kanenobu} T.~Kanenobu, \em{Relations between the
	Jones and Q polynomials of 2-bridge and 3-braid links},
	Math. Ann. {\bf 285} (1989), 115--124.
\bibitem[KS]{KS} M.\ Kidwell and A.~Stoimenow, \em{%
	Examples Relating to the Crossing Number, Writhe,
	and Maximal Bridge Length of Knot Diagrams}, preprint.
\bibitem[Ko]{Kobayashi} M. Kobayashi and T. Kobayashi, {\em On canonical
	genus and free genus of knot}, J. Knot Theory Ramif.\
	{\bf 5(1)} (1996), 77--85.
\bibitem[LM]{LickMil} W. B. R. Lickorish and K.~C.~Millett, {\em A
	polynomial invariant for oriented links}, Topology
	{\bf 26 (1)} (1987), 107--141.
\bibitem[MT]{MenThis} W.~W.~Menasco and M.~B.~Thistlethwaite, \em{%
	The Tait flyping conjecture}, Bull. Amer. Math. Soc. {\bf 
	25 (2)} (1991), 403--412.
\bibitem[Mo]{Morton} H.~R.~Morton, \em{Seifert circles and knot
  	polynomials}, Proc. Camb. Phil. Soc. {\bf 99} (1986), 107--109.
\bibitem[Mo2]{Morton2} \bysame, {\em The Burau matrix and Fiedler's
	invariant for a closed braid}, Topol.\ Appl.
	{\bf 95(3)} (1999), 251--256.
\bibitem[Mo3]{Morton3} \bysame, {\em Threading knot diagrams}, 
	Math. Proc. Cambridge Philos. Soc. {\bf 99(2)} (1986), 247--260.
\bibitem[Mo4]{Morton4} \bysame\ (ed.), {\em Problems},  in
	``Braids'', Santa Cruz, 1986 (J.~S.~Birman and A.~L.~Libgober,
	eds.), Contemp. Math. {\bf 78}, 557--574.
\bibitem[Mo5]{Morton5} \bysame\ (ed.), {\em Problems},  in
	the Proceedings of the International Conference on
	Knot Theory ``Knots in Hellas, 98'', Series
	on Knots and Everything {\bf 24}, World Scientific, 2000.
\bibitem[ME]{MortonElrifai} \bysame\ and E.\ A.\ El-Rifai, \em{%
	Algorithms for positive braids},
	Quart. J. Math. Oxford Ser. 2, {\bf 45(180)} (1994), 479--497.
\bibitem[MS]{MorSho1}  \bysame\ and H. B. Short, \em{Calculating the
	$2$-variable polynomial for knots presented as closed braids},
	J. Algorithms {\bf 11(1)} (1990), 117--131.
\bibitem[MS2]{MorSho2} \bysame\ and \bysame, \em{The
	$2$-variable polynomial of cable knots}, Math. Proc.
	Cambridge Philos. Soc. {\bf 101(2)} (1987), 267--278.
\bibitem[M]{Murakami} J.\ Murakami, \em{The Kauffman polynomial of
	links and representation theory},
	Osaka J. Math. {\bf 24(4)} (1987), 745--758.
\bibitem[Mu]{Murasugi} K.~Murasugi, \em{On the braid index of alternating
	links}, Trans. Amer. Math. Soc. {\bf 326 (1)} (1991), 237--260.
\bibitem[Mu2]{Murasugi2} \bysame, \em{Jones polynomial and classical
        conjectures in knot theory}, Topology {\bf 26} (1987), 187--194.
\bibitem[Mu3]{Murasugi3} \bysame, \em{On closed 3-braids},
	Memoirs AMS {\bf 151} (1974), AMS, Providence.
\bibitem[MP]{MurPrz} \bysame\ and J.~Przytycki, \em{The skein polynomial
        of a planar star product of two links}, Math. Proc. Cambridge
        Philos. Soc. {\bf 106(2)} (1989), 273--276.
\bibitem[N]{Nakamura} T.~Nakamura, \em{Positive alternating links
	are positively alternating}, J. Knot Theory Ramifications
	{\bf 9(1)} (2000), 107--112.
\bibitem[Ro]{Rolfsen} D.~Rolfsen, {\em Knots and links}, Publish
	or Perish, 1976.
\bibitem[Ru]{Rudolph} L.\ Rudolph, \em{Nontrivial positive braids
	have positive signature}, Topology {\bf 21(3)} (1982), 325--327.
\bibitem[Ru2]{Rudolph2} \bysame, {\em Braided surfaces and
	Seifert ribbons for closed braids}, Comment. Math. Helv.
	{\bf 58} (1983), 1--37.
\bibitem[Sl]{Sloane} N.~J.~A.~Sloane,
	{\em The On-Line Encyclopedia of Integer Sequences},
	accessible on the Internet address 
	\verb|http://www.|\lz|research.|\lz|att.com/|\lz|~njas/|\lz|sequences|.
\bibitem[S]{Stallings} J.\ R.\ Stallings, \em{Constructions of fibred
	knots and links}, ``Algebraic and geometric topology'' (Proc.
	Sympos. Pure Math., Stanford Univ., Stanford, Calif., 1976),
	Part {\bf 2}, 55--60.
\bibitem[St]{pos} A.~Stoimenow, {\em Positive knots, closed
	braids and the Jones polynomial}, preprint {\tt math/9805078}.
\bibitem[St2]{gen2} \bysame, {\em Knots of genus two}, preprint.
\bibitem[St3]{fibo} \bysame, \em{Generating functions, Fibonacci
	numbers, and rational knots}, preprint.
\bibitem[St4]{gen1} \bysame, {\em Knots of genus one}, Proc.~Amer.\ 
	Math.~Soc. {\bf 129(7)} (2001), 2141--2156.
\bibitem[St5]{deg} \bysame, {\em Some inequalities between knot
	invariants}, preprint.
\bibitem[Vo]{Vogel} P.~Vogel, {\em Representation of links by braids:
	A new algorithm}, Comment. Math. Helv. {\bf 65} (1990),
	104--113.
\end{thebibliography}
\end{document}